\documentclass[12pt]{amsart}
\usepackage{amscd,amsmath,amssymb,amsfonts}
\theoremstyle{plain}
\newtheorem{thm}{Theorem}
\newtheorem{lem}[thm]{Lemma}
\newtheorem{cor}[thm]{Corollary}
\newtheorem{prop}[thm]{Proposition}
\theoremstyle{definition}
\newtheorem{remark}[thm]{Remark}
\newtheorem{defn}[thm]{Definition}

\numberwithin{thm}{section}
\numberwithin{equation}{section}
\newcommand{\eq}[2]{\begin{equation}\label{#1}#2 \end{equation}}
\newcommand{\ml}[2]{\begin{multline}\label{#1}#2 \end{multline}}
\newcommand{\ga}[2]{\begin{gather}\label{#1}#2 \end{gather}}

\newcommand{\Pic}{{\rm Pic}}

\newcommand{\Spec}{{\rm Spec \,}}

\newcommand{\tr}{{\rm Tr}}

\newcommand{\sD}{{\mathcal D}}
\newcommand{\sE}{{\mathcal E}}

\newcommand{\sI}{{\mathcal I}}

\newcommand{\sK}{{\mathcal K}}

\newcommand{\sO}{{\mathcal O}}

\newcommand{\sV}{{\mathcal V}}
\newcommand{\sW}{{\mathcal W}}

\newcommand{\C}{{\mathbb C}}

\newcommand{\F}{{\mathbb F}}

\renewcommand{\H}{{\mathbb H}}

\renewcommand{\P}{{\mathbb P}}
\newcommand{\Q}{{\mathbb Q}}

\newcommand{\Z}{{\mathbb Z}}
\newcommand{\res}{{\rm res \ }}

\begin{document}

\title{A Formula for Gau\ss-Manin
Determinants} 
\author{Spencer Bloch}
\address{Dept. of Mathematics,
University of Chicago,
Chicago, IL 60637,
USA}
\email{bloch@math.uchicago.edu}

\author{H\'el\`ene Esnault}
\address{Mathematik,
Universit\"at Essen,
Essen,
Germany}
\email{esnault@uni-essen.de}
\date{April 12, 2000}
\thanks{ This work has been partly supported by the NSF grant
DMS-9423007-A1, the DFG Forschergruppe
''Arithmetik und Geometrie'', the Humboldt Foundation and the
University Paris VII}

\begin{abstract}
We give an explicit formula 
for the determinant of the Gau{\ss}-Manin connection for an irregular
connection on a Zariski open set of the projective line $\P^1_K$ 
over a function field $K$
over a field $k$ of  characteristic zero.

\end{abstract}
\maketitle
\begin{quote}
I hasten to write down in verse what I saw then,\\
For the scene lost to sight can't be revived again.\\
\\ \\
\hfill{Su Shi} 
\end{quote}

\section{Introduction}
 
Let $K$ be a function field over a field $k$ of characteristic
0, and let $j: U \subset \P^1_K$ be a Zariski open set of the
projective line. 
We consider a flat connection $(E, \nabla)$   on $U$.
The de Rham cohomology groups $H^i_{DR}(U/K,
\nabla_{/K})$ carry a $K/k$ connection, the Gau{\ss}-Manin connection,
and taking the alternate tensor of the determinant connections
$$\otimes (\det H^i_{DR}(U/K, \nabla_{/K}), {\rm
Gau{\ss}-Manin})^{(-1)^i}, 
$$ one defines the Gau{\ss}-Manin determinant connection,
denoted by $$\det H_{DR}(U/K, \nabla_{/K}).$$ 
This invariant is living in the group of isomorphism classes of $K$-lines
endowed with a connection, which is the abelian group
$$ \Omega^1_K/d\log K^{\times}.$$ 
The aim of this article is to give an explicit formula  for it (see theorem
\ref{mainthm} for a vague formulation, and  \ref{thm:mainthm}
for a  precise one) under a genericity 
assumption on $(E,\nabla)$. Special examples are contained in 
\cite{BE3}.

We comment briefly on the meaning and interest in such a formula. There
is a deep analogy between connections on curves over function fields in
characteristic $0$ and $\ell$-adic sheaves $\sE_\ell$ on curves $U$ over
finite fields $\F_q$. Irregular singular points for the connection
correspond to wild ramification at $\infty$ for the $\ell$-adic sheaf,
and the Gau{\ss}-Manin determinant connection corresponds to the global
epsilon factor
$$\epsilon(\sE_\ell) := \det(-f\ | \ \det R\Gamma(U,\sE_\ell)).
$$
Our hope is that the local formula we obtain for higher rank irregular
connections will suggest local formulas for $\epsilon$-factors extending
the abelian Tate theory.  

The Gau\ss-Manin construction is fairly standard and we do not recall it
in detail. By way of example, we cite two classical formulas (Gau\ss \
hypergeometric and Bessel functions, respectively) :
\begin{gather*}
\frac{\Gamma(b)\Gamma(c-b)}{\Gamma(c)}F(a,b;c;z)
= \int_0^1u^{b-1}(1-u)^{c-b-1}(1-uz)^{-a}du \\
J_n(z) = \frac{1}{2\pi i} \int_{S_0}u^{-n} \exp\
\frac{z}{2}\Big(u+\frac{1}{u}\Big)\frac{du}{u}\quad\quad (S_0 =
\text{circle about
$0$)}.
\end{gather*}

In both cases, the integrand is a product of a solution of a rather
simple degree $1$ differential equation in $u$, the solution being either 
$$u^{b-1}(1-u)^{c-b-1}(1-uz)^{-a}\quad \text{  or}\quad  
u^{-n}\exp\ \frac{z}{2}\Big(u+\frac{1}{u}\Big),
$$
with an algebraic $1$-form ($du$ or $\frac{du}{u}$). The integral is
taken over a chain in the $u$-plane. The resulting functions $F(a,b;c;z)$
and
$J_n(z)$ satisfy  Gau\ss-Manin equations, which are much more
interesting degree $2$ equations in $z$. 

It is not our purpose to go further into the classical theory, but, to
understand the role of the determinant, we remark that in each of the
above cases, there is a second path and a second algebraic $1$-form such
that the two integrals, say $f_1(z)$ and $f_2(z)$, satisfy the same second
order equation. The Wronskian determinant
$$\begin{vmatrix}f_1 & f_2 \\
\frac{df_1}{dz} & \frac{df_2}{dz}\end{vmatrix}
$$
satisfies the degree $1$ equation given by the determinant of
Gau\ss-Manin. It seems to us to be possible using the theory of Stokes
structures to formulate a theory of period integrals for irregular
connections in such a way that the Gau\ss-Manin determinant connection
has as solution the determinant of the period matrix. We hope to return
to this in a future paper. 

Let $X$ be a complete, smooth curve over $K$. For purposes of this article,
we define the group of relative algebraic differential characters
\eq{1.1}{AD^2(X/K):= \H^2(X, \sK_2 \stackrel{d\log}{\longrightarrow}
\Omega^2_X/(\sO_X\otimes\Omega^2_K)) 
}
(The notation here differs from
\cite{E}, as one has factored out $2$-forms coming from the base
and in particular truncated the differential forms of degree
3). The transfer   
\eq{1.2}{f_*: AD^2(X/K) \to AD^1(K)=\Omega^1_K/d\log K^{\times}
} 
maps the group of relative
algebraic differential characters of degree 2 on $X$ to the group of 
algebraic characters of degree 1 on $K$, which is 
the group of connections on $K$. Indeed this is an isomorphism 
(lemma \ref{lem:tr}). But to write connections on $K$ as coming
from differential characters on $X$ allows to single out two
types of classes, global decomposable classes and local
classes, which we discuss in the sequel.

Let $\sD = \sum m_i x_i$ be an effective divisor on $X$, and let $D=\sum
x_i$ be the corresponding reduced divisor. We define a sheaf of
meromorphic $1$-forms 
\eq{1.3}{\Omega^1_X(\sD - D)\subset \Omega^1_X\{\sD\}\subset
\Omega^1_X(\sD) }
as follows. If $z$ is a local parameter at a point $x$ of multiplicity
$m$ in $\sD$, a $1$-form is a section of $\Omega^1_X\{\sD\}$ if it can be
written in local coordinates in the form
\eq{1.4}{\frac{fdz}{z^m} + \frac{\eta}{z^{m-1}}
}
where $f\in \sO_{X,x}$ and $\eta \in \sO_{X,x}\otimes\Omega^1_K$ are
regular at $x$. We define $\Omega^p_X\{\sD\} =\Omega^1_X\{\sD\}\wedge
\Omega^{p-1}_X \subset \Omega^p_X(\sD)$.  There is an exact sequence
\eq{1.5}{0 \to \sO_X(\sD-D)\otimes\Omega^1_K \to \Omega^1_X\{\sD\} \to
\omega_{X/K}(\sD) \to 0 
}
where we write $\omega_{X/K}$ for the sheaf of relative $1$-forms. 

The graded  algebra $\oplus_n \wedge^n(\Omega^1_X\{\sD\})$ is closed
under exterior $d$. Further, writing $\sI = \sO(-\sD) \subset
\sO_X$ for the ideal sheaf, we have
\eq{wedge}{ d\log((1+\sI)^\times )\wedge \Omega^p_X\{\sD\} \subset
\Omega^{p+1}_X. 
}

Let $E$ be a vector bundle on $X$, and let 
\eq{1.6}{\nabla : E \to E\otimes\Omega^1_X\{\sD\}
}
be an absolute connection (i.e. parameters from $K$ are also
differentiated). 
\begin{defn}\label{def1.1}The connection $\nabla$ is {\it vertical}, if the
curvature
$$\nabla^2 : E \to E(*D)\otimes \Omega^2_K.
$$
\end{defn}

We assume our connection $\nabla$ is vertical. (Of course the most
important case is that of a flat connection $\nabla^2=0$.) We write
$\nabla_{/K}:E \to \omega_{X/K}(\sD)$ for the corresponding relative
connection. Viewing
$\sD$ as a nilpotent subscheme of $X$, it is easy to check that $\nabla$
induces a function linear ``polar part'' map 
\ga{1.7}{\nabla_\sD : E|_\sD \to
E\otimes\Big(\Omega^1_X\{\sD\}\big/\Omega^1_X\Big) \qquad\qquad\qquad\ \
\text{(Absolute)}\\ \nabla_{\sD/K} :E|_\sD \to
E\otimes\omega_{\sD/K} = E\otimes(\omega_{X/K}(\sD)/\omega_{X/K})\quad
\text{(Relative)}\notag }
\begin{defn}\label{def1.2} The connection \eqref{1.6} is said to be
{\it admissible} if the relative polar parts map $\nabla_{\sD/K} :E|_\sD \to
E\otimes\omega_{\sD/K}$ is an isomorphism in a singular point of
multiplicity $\ge 2$, and if in a singular point of multiplicity 1,
Deligne's condition \cite{De} that the eigenvalues of the
residue do not belong to $\{0, 1, 2 \ldots\}$ is fulfilled. 
\end{defn}
Notice that the notion of admissibility depends on the extension of $E$
to all of $X$. Although the definition makes reference only to the
relative connection, admissibility depends also on the absolute
connection, which is required to take values in
$\Omega^1_X\{\sD\}\subset \Omega^1_X(\sD)$. It is a local formal
property. 
The motivation for this definition comes from a struture
theorem (\cite{L}, \cite{W}, \cite{Mal}, compare also with \cite{Man},
p.124)
asserting that locally
formally, after ramification of the curve, a flat connection
$(E, \nabla)$ becomes a direct sum of summands $L\otimes
\Lambda$, where $\Lambda$ is a higher rank flat connection with
logarithmic singularities, and $L$ is either a rank 1 trivial
connection or a rank 1 flat connection with
multiplicity $\ge 2$ (strictly speaking, this is proven in
\cite{Mal} only over $\C$, outside of a Baire set). Since rank 1
vertical connections are admissible (\cite{BE2}, lemma 3.1), we
see that up to ramification, any flat connection is locally of
sum of admissible and logarithmic connections. 

For a rank 1 vertical connection we may view
\eq{1.9}{\nabla_{\sD/K}: E|_\sD \cong E|_\sD\otimes
\omega_{X/K}(\sD)|_\sD  
}
as defining a trivialization of $\omega_{X/K}(\sD)|_\sD$, i.e. a class
$(\omega_{X/K}(\sD), \nabla_{\sD/K})\in \Pic(X,\sD)$. There is a
cohomological pairing ($j:X-D \hookrightarrow X$)
\begin{gather} \label{eq1.10}
\{\ \,,\ \}:\Pic(X,\sD)\otimes \H^1(X,j_*\sO_{X-D}^\times \to
\Omega^1_X\{\sD\}) \\
\to AD^2(X)=\H^2(X, \sK_2\xrightarrow{d\log}\Omega^2_X) \notag
\end{gather}
and in the rank $1$ case (see main theorem of \cite{BE2})
\begin{gather}\label{eq1.11}
\det H_{DR}(U/K, \nabla_{/K}) = -
f_*\{(\omega_{X/K}(\sD),\nabla_{\sD/K}),(E,\nabla)\}\\
 \in
\Omega^1_K/d\log K^{\times} \otimes \Q.\notag
\end{gather}
In higher rank, however, we have examples of connections $(E,\nabla)$
with $(\det(E),\det(\nabla))$ trivial but non-trivial Gau\ss-Manin
determinant connection (see \cite{BE3},
remark 3.3, equation 3.27, and remark \ref{rmk2.12}). 

Still, the choice of a meromorphic section $s\in
\omega_{X/K}(\sD)\otimes K(X)$ which generates the sheaf at the
points of $\sD$ defines a rigidification $c_1(\omega_{X/K}(\sD),
s) \in \Pic(X, \sD)$, and allows to define a class 
$$\{c_1(\omega_{X/K}(\sD), s), \det(E,\nabla)\}\in AD^2(X).$$
We refer to it as the global factor (see \eqref{eq2.5}).

In each singularity, we define local factors (see proposition
\ref{prop:ind}) which play the r\^{o}le of the local epsilon
factors defined to express the global epsilon of an $\ell$-adic
sheaf. If $A_i= g_i\frac{dz}{z^m} + \frac{\eta_i}{z^{m-1}}$
is the local equation of $\nabla$ in a local basis around the
point $a_i\in D$, then one defines $\tr dg_i g_i^{-1}A_i \in
AD^2(X/K)$. We also define a 2-torsion local factor $\frac{1}{2}
d\log (\det g_i(a_i))) \in AD^2(X/K)$ (see definition \ref{defn:tau}).

The main theorem of this article says
\begin{thm}\label{mainthm} Let $(E,\nabla)$ be an admissible
connection on $\P^1_K$ having  at least one point of multiplicity $\ge
2$. Then 
\begin{gather*}
\det H_{DR}(U, \nabla_{/K})= -f_*\{c_1(\omega_{X/K}(\sD), s), \det(E,\nabla)\}
 +\\
\sum_i \res\!\!_{a_i}\, \tr\Big(dg_ig_i^{-1}A_i\Big) +\frac{1}{2}\sum_i
m_i d\log(\det(g_i(a_i)))\\
\in \Omega^1_K/d\log(K^\times).
\end{gather*}
\end{thm}
(See theorem \ref{thm:mainthm} for  a slightly more precise 
formulation). 

The cases rank $(E)=1$, resp. $\nabla$ with logarithmic poles,
were considered in  the earlier articles \cite{BE2}, resp.
\cite{BE1}, except that for the rank 1 case, our results  did not
include torsion.  In those two cases,  there is a well defined class
$\gamma(E,\nabla) \in AD^2(X)$ such that $f_*\gamma(E,\nabla)=\det
H_{DR}(U/K, \nabla_{/K})$.  In the higher rank, non-logarithmic
case considered here, one still has the global factor in $AD^2(X)$, but
the local factors are well defined in $AD^2(X/K)$ only (see
proposition \ref{prop:ind}). 

It would be of great interest to find some variant of this formula which
applied to epsilon factors for $\ell$-adic sheaves.

We now discuss the proof of the main theorem.
Let $K$ be a field in characteristic 0. A classical theorem by Euler 
(\cite{Se}, III, 6, lemma 2) asserts that if $g \in K[t]$ is a
polynomial, and $h \in L=K[u]/gK[u]$, then the trace of 
multiplication by $h$, viewed as a $K$-linear map from $L$ to itself,
is computed by 
$$-\res\!\!_{u=\infty}dgg^{-1}h.
$$ 
We define a
generalization of this in the non-commutative situation as follows. Let
$V$ be a finite dimensional $K$-vectorspace, and $g=\sum_{i=0}g_iu^i \in
{\rm
  End}(V)[u]$ be a polynomial with coefficients in the endomorphisms
of $E$, such that the leading coefficient $g_m \in {\rm Aut}(V)$ 
is invertible. The invertibility of $g_m$ allows  us to write an
element of
$W=V[u]/gV[u]$ as the class of an element $\sum_{r=0}^{m-1}v_i
u^i$, where $v_i\in V$, yielding 
a splitting $\sigma: W 
\to V[u]$ of the natural projection $p: V[u] \to W$. For any $h\in 
{\rm End}(V)[u]$, $\phi(h):= p\circ h\circ\sigma  : W \to W$ will have
a trace. Proposition \ref{prop5.1} says
$$ \tr_{V[u]/gV[u]}(\phi(h)) = - \tr_V \res\!\!_{u=\infty}(dgg^{-1}h).
$$

The second point is to relate the trace of this linear operator with
the trace of a differential operator. The crucial case to
understand is that of a connection on a trivial bundle $E=V\otimes_K
\sO_{\P^1}$ on
$\P^1_K$. Such a connection is given by a matrix which  has the shape
\eq{1.12}{
\sum_{i=1}^N \sum_{r=1}^{m_i}
\frac{g_r^{(i)}d(t-a_i)}{(t-a_i)^r} + 
\eta = g + \eta,}
where $g_r^{(i)}\in {\rm End}(V)$ and $\eta \in {\rm End}(V)\otimes
\Omega^1_K\otimes
\sO_{\P^1_K}(*D)$. We write $g$ also for the corresponding matrix of
relative forms, so $\nabla_{/K} = d+g$. We fix a certain
finite-dimensional vector subspace
$\sigma:H \hookrightarrow H^0(\P^1_K,V\otimes_K\omega(*D))$ such that
composition with the natural projections give isomorphisms 
\ga{}{H\stackrel{\sigma}{\to}H^0(\P^1_K,
V\otimes_K\omega(*D)) \stackrel{p_\nabla}{\to}H^1_{DR}(U/K,
\nabla_{/K}) \\
=H^0(\P^1_K,
V\otimes_K\omega(*D))/{\rm Im}\nabla_{/K}\notag \\
H\stackrel{\sigma}{\to}H^0(\P^1_K,
V\otimes_K\omega(*D)) \stackrel{p_g}{\to}H^0(\P^1_K,
V\otimes_K\omega(*D))/{\rm Im} \ g   \notag
}
The operators
\ga{}{\eta_\nabla
:= (p_\nabla\circ\sigma)^{-1}p_\nabla\circ\eta\circ\sigma:H
\to H\otimes\Omega^1_K
\notag \\
\eta_\gamma 
:=(p_g\circ\sigma)^{-1}p_g\circ\eta\circ\sigma:H \to H\otimes\Omega^1_K
\notag
 }
will be referred to as the (Gau\ss-Manin) de Rham operator and Higgs
operator respectively. The traces of these operators play a central role
in the Gau\ss-Manin determinant, and the remarkable fact is that for an
admissible, vertical connection on a trivial bundle on $\P^1_K$ one finds
\eq{1.13}{\tr(\eta_\nabla -\eta_\gamma) \equiv \frac{1}{2}\sum_i m_i
d\log(\det(g_{m_i}^{(i)})) \mod d\log(K^\times).
}
(Here the $g^{(i)}_{m_i}$ are as in \eqref{1.12}.) This result is 
theorem \ref{thm:trop}. It is reminiscent of Hitchin's 
comparison of de Rham and Higgs twisted cohomologies on projective
manifolds, and of Kontsevich's theorem comparing de Rham and
Higgs cohomology of $df$, where $f$ is a regular function on a manifold.
Algebraically, if 
$$(\frac{g_m}{z^m} + \frac{g_{m-1}}{z^{m-1}} + \ldots)dz +
\frac{\eta_{m-1}}{z^{m-1}} + \frac{\eta_{m-2}}{z^{m-2}} + \ldots,$$
represents the polar part of our admissible, vertical connection at a
point $z=0$, the essential result (proposition \ref{prop:iden}) is that
$$\tr \sum_{s=0}^{m-1} g_m^{-1}[g_{m-s}, \eta_s]$$
is identically vanishing.
We must confess that, even  after performing
 the computation, we don't really understand its meaning.

{\it Acknowledgements}: It is a pleasure to thank A. Beilinson,
C. Sabbah and T. Saito for interesting discussions related to
the topics discussed in this article.

\section{ Admissible Connections}
Let $K$ be a function field over a field $k$ of characteristic
0, $f: X\to \Spec K$ be a smooth projective curve, $j: U\subset
X$ a non-trivial Zariski open set such that the closed points
of $D=X\setminus U$ are $K$ rational points, and $(E, \nabla)$ a
global connection of rank $r$ on $X$ which is regular on $U$. Barring
express mention to the contrary, we shall always assume $\nabla$ to be
vertical (definition
\ref{def1.1}). We shall need a small generalization of the notion of
admissibility introduced in definition \ref{def1.2}. 
\begin{defn} \label{defn:psadm}
The connection
$(E, \nabla)$ is {\it pseudo-logarithmic }at  $x \in X\setminus
U$ if the local equation of $\nabla$ in some basis of $E$ has the shape
$$ A= g\frac{dz}{z} + \frac{\eta}{z},$$
where $z$ is a local parameter around $x$, $\eta \in M(r\times
r, \Omega^1_K\otimes \sO_X)$, $g \in GL(r,
\sO_{x})$. The connection $(E, \nabla)$ is {\it pseudo-admissible} if
it is admissible in singularities with multiplicities $m\ge 2$
and pseudo-logarithmic in points with multiplicities $m=1$.
\end{defn}

We will use a very special simple shape of pseudo-logarithmic
singularities, which we single out in the following definition.
\begin{defn} \label{defn:spe}
A {\it special pseudo-logarithmic point} of a connection
$(E,\nabla)$ is pseudo-logarithmic, and there is 
local basis $$(e_\nu)=( (e_1,\ldots, e_s), (e_{s+1},
\ldots, e_r))$$ with respect to which the block-matrix of the
connection has the shape 
\begin{gather*}\Big(
\begin{array}{ll}
A + m\frac{dz}{z} & zB\\
\frac{C}{z} & D + n\frac{dz}{z}
\end{array}\Big),
\end{gather*}
where the connection matrix
\begin{gather*}\Big(
\begin{array}{ll}
A  & B\\
C & D 
\end{array}\Big)
\end{gather*}
has no poles and $m,n \in k$.

A connection is special pseudo-logarithmic if it is admissible,
and special pseudo-logarithmic in pseudo-logarithmic points. 
\end{defn}

Working with pseudo-admissible connections will enable us to reduce the
Gau\ss-Manin determinant computation for a general $(E,\nabla)$ on $\P^1$
to the case where $E\cong \sO_{\P^1}^{\oplus r}$. We use the following:

\begin{thm} \label{thm:red} [Compare with \cite{BE1}, lemma 4.2 and
    reduction 4.1]
Let $(E,\nabla)$ be an admissible connection on $\P^1_K$ having a
singularity of multiplicity $\ge 2$.
Then there are finitely many points $p_i\in U(K)$, such that if
$\lambda: V=U\setminus \{p_i\} \to \P^1_K$ denotes the open embedding, then
$(E|_V, \nabla|_V)$ extends to a special pseudo-admissible connection
$(\oplus_1^r 
\sO_{\P^1_K}, \nabla)$ on $\P^1_K$ such that 
\begin{gather*} 
(\oplus_1^r \sO_{\P^1_K})
\xrightarrow{\nabla_{/K}} \omega(\sD + \sum_i p_i)\otimes  
(\oplus_1^r \sO_{\P^1_K})\\ \to (\lambda_*E_V 
\xrightarrow{\nabla_{/K}}\lambda_*( \omega\otimes 
E_{V}))
\end{gather*}
 is a quasiisomorphism.
\end{thm}
\begin{proof} 
Without loss of generality, we may assume 
that $\infty$ is a smooth point of the
connection. Let $x$ be a point of multiplicity $\ge 2$, and let $z$ be a
local coordinate at $x$. Note the effect of twisting (i.e. replacing
$E$ by
$E(Nx)$) is to replace the local connection matrix $A$
at $x$ with respect to a basis $e_i$ with $A-N\frac{dz}{z}I$ for the basis
$\frac{e_i}{z^N}$. In particular, $x$ will remain an admissible
singularity for the new connection. After such a twist, we may assume $E=
\oplus_{i=1}^r
\sO(n_i)$, with $0\le n_1 \le n_2 \le \ldots$.  We argue by induction on
$n_r-n_1$. If $n_r-n_1=0$, we replace $E$ by $E(-n_1x)$ and argue as
above.

Assume $n_r-n_1>0$. Let
$E'=\oplus_{i=1}^{r-1}\sO(n_i) \oplus \sO(n_r-1)$, and embed $E'$ in $E$
via $\sO((n_r-1)\infty)\to  \sO(n_r\infty)$.  If $z$ is a local parameter
at $\infty$, and $e_\nu$ is  a local basis of $\sO(n_\nu)$ at $\infty$,
then $((e_\nu, \mu \le r-1), ze_r)$ is a local basis of $E'$, and if
\begin{gather}
\Big(\begin{array}{ll}
A & B\\
C & D
\end{array}\Big)
\end{gather}
 is the local block matrix of the connection $\nabla$ in the basis
$((e_1,\ldots,e_{r-1}), e_r)$, then
\begin{gather}
\Big(\begin{array}{ll}
A & zB\\
\frac{C}{z} & D + \frac{dz}{z}
\end{array}\Big)
\end{gather}
is the local block matrix of the connection in the basis 
$((e_1,\ldots,e_{r-1}), ze_r)$. 
If $C$ has a local  expansion $C= C_0 + C_1 z + \ldots$, then the polar part
of this connection is
\begin{gather}
\Big(\begin{array}{ll}
0 & 0\\
\frac{C_0}{z} &  \frac{dz}{z}
\end{array}\Big).
\end{gather}
Thus replacing now $E'$ by
$E''=E'(2\infty)\cong \oplus_{i=1}^{r-1}\sO(n_i+2)\oplus
\sO(n_r+1)$, the local equation of the connection at $\infty$ 
becomes
\begin{gather}\label{plspe}
\Big(\begin{array}{ll}
A -2\frac{dz}{z} & 0\\
\frac{C}{z} & D -\frac{dz}{z}
\end{array}\Big),
\end{gather}
and therefore, is pseudo-admissible. 
On the other hand, $n_r-n_1$ has decreased. We conclude by induction. 
\end{proof}

Next we describe the class
$\gamma(E,\nabla) \in AD^2(X/K) :=
\H^2(X,\sK_2\stackrel{d\log}{\longrightarrow}\Omega^2_X/\sO_X\otimes\Omega^2_K)$
from theorem \ref{thm:mainthm}. If we choose a meromorphic section $s$ of
$\omega_{X/K}(\sD)$ which generates this sheaf in a neighborhood of
$\sD$, we may view $s$ as defining a trivialization of
$\omega_{X/K}(\sD)|_\sD$, i.e. a class $(\omega(\sD),s) \in \Pic(X,\sD)$.
As in \eqref{eq1.10}, we may consider the product
\begin{gather}\label{eq2.5}
\{(\omega(\sD),s),(\det(E),\det(\nabla))\} \in
AD^2(X). \end{gather}
We refer to this class as the global factor.

Fix a basis $e_i$ for $E$ in a neighborhood of $\sD$. the choice of $e_i$
determines a local connection matrix $A$, so $\nabla = d+A$. Let
$\sO_{X,\sD}$ denote the semi-local ring of functions regular at all
points of $\sD$. The choice of $s$ determines $g\in GL_r(\sO_{X,\sD})$
such that the relative connection $\nabla_{/K} = d +gs$. Note the
hypothesis of pseudo-admissibility insures that $g$ is invertible. 

The basic local invariant we consider is 
\eq{2.6}{\tr(dgg^{-1}A) \in H^0\Big(X,\Omega^2_X\{\sD\}\Big/\big(\Omega^2_X +
\sO_X(\sD-D)\otimes\Omega^2_K\big)\Big). 
}
The boundary map from the exact sequence
\ml{2.7}{ 0 \to \Omega^2_X/\sO_X\otimes\Omega^2_K \to \Omega^2_X\{\sD\}\Big/
\Big(\sO_X(\sD-D)\otimes\Omega^2_K\Big) \\
\to \Omega^2_X\{\sD\}\Big/
\Big(\Omega^2_X +\sO_X(\sD-D)\otimes\Omega^2_K\Big) \to 0
}
together with the evident map $H^1(X, \Omega^2_X/\sO_X\otimes\Omega^2_K) \to
AD^2(X/K)$ enables us to define an element, which we denote by abuse of notation
\eq{2.8}{\tr(dgg^{-1}A) \in AD^2(X/K).
}

\begin{prop}\label{prop:ind} 
Let $(E,\nabla)$ be a pseudo-admissible vertical 
connection on $X$. The element
$\tr(dgg^{-1}A)$ \eqref{2.8} is independent of the choice of local bases around
$\sD$. The element
\eq{2.9}{-\{(\omega(\sD),s),(\det(E),\det(\nabla))\} + \tr(dgg^{-1}A) \in
AD^2(X/K) 
}
is independent both of the choice of local bases and the trivializing
meromorphic section $s$ of $\omega(\sD)$. 
\end{prop}
\begin{proof}We show first that $\tr(dgg^{-1}A)$ is independent of the local
bases. We work locally around a point $x$ which is a singular point of the
connection with multiplicity $m$. We assume first that $s=\frac{dz}{z^m}$ for
a local coordinate, and that the connection matrix is
$$A = \frac{gdz}{z^m} + \frac{\eta}{z^m}.
$$ 
with $g\in GL_r(\sO)$ and $\eta\in  M_r(\sO\otimes\Omega^1_K)$. 
(This includes both the admissible and pseudolog cases.) We take
a gauge transformation of the form $A\mapsto \phi A \phi^{-1}
+d\phi\phi^{-1}$ with $\phi\in GL_r(\sO)$. This amounts to
$$g \mapsto \phi g\phi^{-1} + z^m \frac{d\phi}{dz}\phi^{-1}; \ \eta
\mapsto \phi \eta \phi^{-1} + z^m d_K \phi \phi^{-1}. 
$$  
Here $d=d_z + d_K$. We claim 
$$\tr(dgg^{-1}A)\in
\Omega^2_X(\sD)\Big/\Big(\sO(\sD)\otimes\Omega^2_K +
\Omega^2_X\Big)
$$ 
is invariant. We compute (writing $T(A) := \tr(dgg^{-1}A)$, and 
computing modulo
$\Omega^2_X$)
\begin{multline} T(\phi A\phi^{-1} + d\phi \phi^{-1}) = \\
\tr\Big(d(\phi g\phi^{-1} + z^m \frac{d\phi}{dz}\phi^{-1})(\phi
g\phi^{-1} + z^m \frac{d\phi}{dz}\phi^{-1})^{-1}(\phi A\phi^{-1} +
d\phi\phi^{-1})\Big) \equiv \\
\tr\Big(d(\phi g\phi^{-1} + z^m \frac{d\phi}{dz}\phi^{-1})\phi g^{-1}
\phi^{-1}\phi A \phi^{-1}\Big) \equiv \\
\tr\Big((d\phi \phi^{-1} + \phi dgg^{-1}\phi^{-1} - \phi g\phi^{-1}
d\phi\phi^{-1} \phi g^{-1} \phi^{-1} \\ + mz^{m-1} \frac{d\phi}{dz}
g^{-1}\phi^{-1}dz ) \phi A \phi^{-1}\Big) \equiv \\
\tr\Big(\phi^{-1}d\phi A + dgg^{-1} A -g\phi^{-1}d\phi g^{-1} A +
mz^{m-1}\phi^{-1} \frac{d\phi}{dz} g^{-1} dz\wedge A\Big) \equiv \\
\tr\Big( dgg^{-1} A + \phi^{-1} d\phi (A - g^{-1} A g) + mz^{m-1}
\phi^{-1} \frac{d\phi}{dz} g^{-1} dz\wedge A\Big) \equiv \\
\tr\Big(dgg^{-1} A + \phi^{-1}\frac{d\phi}{dz} g^{-1}dz(gA - Ag +
mz^{m-1} A)\Big) 
\end{multline}
(Note that $gA-Ag$ has entries in $\Omega^1_K\otimes K(X)$, justifying
replacing $d\phi$ by $\frac{d\phi}{dz}dz$.) To show invariance, it
will suffice to show
$$dz\wedge(gA-Ag+mz^{m-1}A)\equiv 0 \mod \Omega^2_X.
$$
This expression can be written
$$(*) = \frac{dz}{z^m}\wedge [g,\eta] + m\frac{dz}{z}\wedge\eta. 
$$
Verticality gives
$$dg\wedge\frac{dz}{z^m} + m\eta\wedge\frac{dz}{z^{m+1}} +
\frac{d_z\eta}{z^m} = \frac{dz}{z^{2m}}[g,\eta]
$$
Multiplying through by $z^m$, the expression $(*)$ above becomes  
$$(*) \equiv m\frac{dz}{z}\wedge\eta + m\eta\wedge\frac{dz}{z} \equiv 0
\mod \Omega^2_X. 
$$

It remains to show independence of $s$. Let $s'=fs$ where $f$ is meromorphic on
$X$ and invertible on $D$. Consider a diagram 
\begin{tiny}
\eq{2.10}{\minCDarrowwidth.5cm \begin{CD} H^0(\sO_\sD^\times) @>\partial >>
 H^1((1+\sI_\sD)^{\times})=\Pic(X,\sD) @>>> H^1(\sO_{X}^\times)=\Pic(X) \\
\times @. \times @. \times  \\
\Omega^1_X\{\sD\}\big/\Big(\Omega^1_X + d\log j_*\sO_{X-D}^\times\Big)
@<<<
\H^1(j_*\sO_{X-D}^\times \to \Omega^1_X\{\sD\}) @<<< \H^1(\sO_X^\times \to
\Omega^1_X) \\
@VVV @VVV @VVV \\
\Omega^2_X\{\sD\}\Big/\Big(\Omega^2_X + d\log j_*\sK_{2,\sO_{X-D}}\Big) @>>>
AD^2(X/K) @= AD^2(X/K) 
\end{CD}
}
\end{tiny}
The classes of the rigidified bundles $(\omega(\sD),s)$ and $(\omega(\sD),fs)$
differ by $\partial(\tilde{f})$, where $\tilde{f}\in H^0(\sO_\sD)$ is the image
of $f$. It follows that (with notation as above)
\eq{2.11}{ \{(\omega(\sD),fs)(\omega(\sD),s)^{-1},(\det(E),\det(\nabla))\} =
\tr(\frac{d\tilde{f}}{\tilde{f}}A) }
Replacing $s$ with $fs$ in \eqref{2.9}, this gives the desired invariance. 
\end{proof}

To complete the construction of the class $\gamma(E,\nabla)\in
AD^2(X/K)$ we need 
one more invariant. Let $x\in D$ and assume $m \ge 2$, where $m$ is the
multiplicity of $x$ in $\sD$. Let $z$ be a local coordinate at $x$. Write the
local relative connection $\nabla_{/K} = d+g\frac{dz}{z^m}=
\frac{g_mdz}{z^m}+\frac{g_{m-1}dz}{z^{m-1}}+\ldots$ with $g_m\in GL_r(K)$. 
\begin{defn} \label{defn:tau}If the multiplicity is $\ge 2$, the invariant
$$\tau_x(E,\nabla):= \frac{m}{2}d\log(\det(g_m))\in
\frac{1}{2}d\log(K^\times)\big/d\log(K^\times)$$
is associated to the rank 1 quadratic form $$\det 
g_m^{m} \in K^\times/(K^\times)^2
\cong H^1(K, \Z/2\Z).$$ In a point of multiplicity 1, we set 
$\tau_x(E, \nabla)=1$.
\end{defn}

A
change of local gauge replaces $g$ with $hgh^{-1}+z^m\frac{dh}{dz}h^{-1}$ with
$h\in GL_r(K[[z]])$. It follows that $\det(g_m)$ is invariant under gauge
transformation of $m\ge 2$. 
On the other hand, replacing $z$ with $z' = uz$ leads to $g_m'
= u^{m-1}g_m$, whence $\det(g_m') = u^{r(m-1)}\det(g_m)$ and
$$\frac{m}{2}d\log(\det(g_m')) = \frac{m}{2}d\log(\det(g_m)) +
\frac{rm(m-1)}{2}d\log(u),
$$
so the definition is independent of the choice of $z$. 

For $x\in X(K)$, we define a map 
\eq{2.12}{\rho_x : K^\times/K^{\times 2} \to AD^2(X/K)
}
as follows. One has a map of
exact sequence of complexes
\eq{2.13}{\minCDarrowwidth.5cm\begin{CD}\sK_{2 X} @>>> j_{x*}\sK_{2,X-\{x\}}
@>>> i_{x*}K^\times \\
@VVV @VVV @VVaV \\
\Omega^2_X/\sO_X\otimes\Omega^2_K @>>> \Omega^2_X\{x\}/\sO_X\otimes\Omega^2_K
@>>> \Omega^2_X\{x\}/\Big(\sO_X\otimes\Omega^2_K + \Omega^2_X\Big).
\end{CD}
}
Write $z$ for a local coordinate and let $a$ be as in \eqref{2.13}. The
mapping 
\eq{2.14}{K^\times \otimes \Q/\Z \to \text{coker}(a);\quad \kappa\otimes
\frac{1}{n} \mapsto \frac{1}{n}\frac{dz}{z}\wedge \frac{d\kappa}{\kappa} 
}
is well-defined. We compose this with the boundary map from \eqref{2.13} to
define $\rho_x$. 
\begin{defn} \label{defn:gam}
Let $(E,\nabla)$ be a pseudo-admissible connection as above. Then
\ml{2.16}{\gamma(E,\nabla) = -\{(\omega(\sD),s),(\det(E),\det(\nabla))\} \\
+ \tr(dgg^{-1}A) +\sum_{x\in \sD, m_x\ge 2}\rho_x(\tau_x(E,\nabla))\in AD^2(X/K). 
}
\end{defn}
We continue to assume $f:X \to \Spec(K)$ is a smooth, projective curve. The
transfer map $f_* : AD^2(X/K) \to \Omega^1_K/d\log(K^\times)$ is defined as
follows. We remark that
$H^2(X,\sK_2)=(0)$ and $\Omega^2_X/\sO_X\otimes\Omega^2_K \cong
\omega_{X/K}\otimes\Omega^1_K$, and we define $f_*$ from the diagram
\eq{2.17}{\minCDarrowwidth.5cm\begin{CD} H^1(X,\sK_2) @>>>
H^1(X,\omega_{X/K})\otimes \Omega^1_K @>>> AD^2(X/K) @>>> 0 \\
@VV\tr V @VV\cong V @VV f_* V \\
K^\times @>>> \Omega^1_K @>>> \Omega^1_K / d\log(K^\times) @>>> 0
\end{CD}
}
The check that with $\gamma(E,\nabla)$ defined as in \eqref{2.16},
$f_*\gamma(E,\nabla)$ has the form as in theorem \ref{mainthm}
 is straightforward and will be
omitted. 
Since the trace map $\tr: H^1(X, \sK_2) \to K^{\times}$, which
is simply defined on the generators $\oplus_{x}\lambda_x\in 
\oplus_{x\in X^{(1)}}
K(x)^{\times}$ by $\tr_{[K(x):K]}d\log \lambda_x$, is surjective, we
obtain the 
\begin{lem} \label{lem:tr}
The transfer map $$f_*: AD^2(X/K)\to AD^1(K)=\Omega^1_K/d\log
K^\times$$ is an isomorphism.
\end{lem}

Now we are in the position to give a slightly more
precise formulation of our
main theorem.
\begin{thm} \label{thm:mainthm}
Let $(E,\nabla)$ be a special pseudo-admissible  connection on $\P^1_K$,
smooth over $\emptyset \neq U\subset \P^1_K$,
such that the singularities $\P^1_K\setminus U=D$
of $\nabla$ consist of  $K$-rational points.
Then, with the notation of definition \ref{defn:gam}, one has
$$\det H_{DR}(U/K, \nabla_{/K})=f_*\gamma(E,\nabla) \in
AD^1(K)=\Omega^1_K/d\log K^\times.
$$
\end{thm}

Finally, we take a moment to point out some simple consequences of theorem
\ref{thm:mainthm}.

\begin{remark} \label{rmk2.12}(i).  
Let $(E,\nabla)$ be an admissible connection
on
$\P^1_K$, and let $(E^\vee,\nabla^\vee)$ be the dual connection. then
\eq{}{\det H_{DR}(E,\nabla) \cong  \det H_{DR}(E^\vee, \nabla^\vee)^{(-1)}.
}  
Indeed, replacing $E$ with $E^\vee$ replaces $g$ with $-{}^tg$ and $A$
with $-{}^tA$. Verticality and admissibility imply that $[g,A]$ has no
pole, so 
\begin{multline*}\tr(d(-{}^tg)(-{}^tg)^{-1} (-{}^tA)) =
-\tr({}^t(g^{-1}dg){}^tA) = -\tr(g^{-1}dgA) \\
= -\tr(dgAg^{-1}) = -\tr(dgg^{-1}A). 
\end{multline*}
 (ii). In
\cite{BE3}, remark 3.3, equation 3.27, we give an example of an admissible
connection on
$\P^1_K$ with trivial determinant, for which the determinant of the
Gau{\ss}-Manin connection is not torsion. More generally, given the main theorem
of \cite{BE2}, if $(L, \nabla_L)$ and $(M,\nabla_M)$ are two
connections, say on $\P^1_K$, with the same singularities $\sD$,
which are generic enough 
so that the singularites of $(L\otimes M, \nabla_L\otimes \nabla_M)$
are exactly $\sD$ as well, then $\det A$, with 
\begin{gather}
A=\Big( \begin{array}{lll}
L\otimes M & 0 & 0 \\
0 & L^{-1} & 0 \\
0 & 0 & M^{-1}
\end{array}\Big)
\end{gather}
is trivial, while the main theorem of \cite{BE2} says that
the Gau{\ss}-Manin determinant is computed by
\begin{gather}
c_1(\omega(\sD), \Gamma_L + \Gamma_M)\cdot c_1(\nabla_L + \nabla_M)\\
-c_1( \omega(\sD), \Gamma_L)\cdot c_1(\nabla_L) -
c_1( \omega(\sD), \Gamma_M)\cdot c_1(\nabla_M),\notag
\end{gather}
where $\Gamma_L$ and $\Gamma_M$ are the principal parts of
$\nabla_L$ and $\nabla_M$. For generic $L$ and $M$, this won't
vanish. 
\end{remark}

\section{Higgs and de Rham Traces}\label{sect3}

In this section we introduce the concepts of 
Higgs and de Rham operators associated to a vertical
pseudo-admissible connection on the trivial bundle on $\P^1_K$,
and analyse the difference between the traces of those
operators. 

It will  be convenient to write $E=V\otimes \sO$ with
$V=\Gamma(\P^1,E)$. Let $(e_\mu), \mu=1,\ldots,r $ be a basis of
$V$. The connection $\nabla$ will be a vertical pseudo-admissible
connection on $E$ with poles on $D=\sum_{i=1}^N (a_i)$ where for
simplicity we take $a_i\in \P^1(K)$. We assume $\infty\not\in D$. The
relative connection is given by
\begin{equation}\label{3.1}\nabla_{/K} e_\mu = \sum_{i=1}^N \sum_{r=1}^{m_i}
\frac{g_r^{(i)}(e_\mu)dt}{(t-a_i)^r}. 
\end{equation}
Here the $g^{(i)}_r$ are matrices with entries in $K$. The
admissibility condition implies that $g^{(i)}_{m_i}$ is
invertible over $\sO_{\P^1_K, a_i}$. Regularity of the
connection at infinity means 
\begin{gather}
\sum_{i=1}^N g^{(i)}_1 = 0.
\end{gather}

\begin{defn} \label{defn:M_i}
Let $M_i=m_i-1$ if $m_i\ge 2$, and $M_i=m_i=1$ else.
\end{defn}
The absolute connection has the following equation  
\begin{equation}\label{3.2}\nabla e_\mu = \sum_{i=1}^N \sum_{r=1}^{m_i}
\frac{g_r^{(i)}(e_\mu)d(t-a_i)}{(t-a_i)^r} + \sum_{i=1}^N
\sum_{r=1}^{M_i} \frac{\eta^{(i)}_r(e_\mu)}{(t-a_i)^r} + \eta_0(e_\mu),
\end{equation}
where the $\eta$ are matrices with  entries in $\Omega^1_K$. 

\begin{defn} \label{defn:gammaK}
We define $\gamma_K$ by $\nabla_{/K}=d+\gamma_K$, thus concretely
$$\gamma_K= \sum_{i=1}^N \sum_{r=1}^{m_i}
\frac{g_r^{(i)}dt}{(t-a_i)^r}.$$
We define $\eta$ by
$$\eta=\sum_{i=1}^N\sum_{r=1}^{M_i} \frac{\eta^{(i)}_r}{(t-a_i)^r} + \eta_0,$$
and  $\gamma$ by
$$\gamma=\sum_{i=1}^N \sum_{r=1}^{m_i}
\frac{g_r^{(i)}d(t-a_i)}{(t-a_i)^r},$$
so that $\nabla = d+\gamma+\eta$. 
\end{defn}
We have natural identifications
\begin{gather} \Gamma(E(*D)) =
V[\frac{1}{t-a_1},\dotsc,\frac{1}{t-a_N}] \\
\label{strict}
\Gamma(E\otimes\omega(*D))\varsubsetneq 
V[\frac{1}{t-a_1},\dotsc,\frac{1}{t-a_N}]dt
\end{gather} 
where the strict inclusion in \eqref{strict} comes from the requirement of
no poles at infinity. 

The following lemma will be useful 
in the sequel.

\begin{lem}\label{lem3.3} For integers $r,s\ge 1$ one has a formal identity
$$\frac{1}{(t-a)^r(t-b)^s} = \sum_{p=1}^r\frac{A_p(a,b)}{(t-a)^p} +
\sum_{q=1}^s\frac{B_q(a,b)}{(t-b)^q}. 
$$
One has
$$A_r = (a-b)^{-s};\quad B_s = (-1)^r(a-b)^{-r} = (b-a)^{-r}.
$$
The partial fraction expansion of $\frac{1}{(t-b)^s(t-a)}$ begins
$$\frac{1}{(t-b)^s(t-a)} = \frac{(a-b)^{-s}}{t-a} - 
\frac{(a-b)^{-s}}{t-b}+\ldots  
$$
In particular, we have
\begin{multline}\frac{1}{(t-b)^s}\Big(\frac{d(a-b)\wedge dt}{t-a} - 
\frac{d(a_N-b)\wedge dt}{t-a_N}\Big) \\
 =  (a-b)^{-s}d(a-b)\wedge dt\Big
 ( \frac{1}{t-a}-\frac{1}{t-a_N}\Big) \\  
+ \Big((a_N-b)^{-s}d(a_N-b)\wedge dt -
 (a-b)^{-s}d(a-b)\wedge dt\Big) \\
\times \Big(\frac{1}{t-b}-
\frac{1}{t-a_N}\Big) + \text{   terms involving $\frac{1}{(t-b)^r}$ for $r\ge
2$} .  
\end{multline}
\end{lem}
\begin{proof} We have
\begin{multline}\frac{1}{(t-b)^s(t-a)^r} \\
= \frac{1}{(s-1)!(r-1)!}
\Big(\frac{d}{db}\Big)^{s-1}\Big(\frac{d}{da}\Big)^{r-1}
\Big(\frac{(a-b)^{-1}}{t-a} - \frac{(a-b)^{-1}}{t-b}\Big).   
\end{multline}
The formulas in the lemma follow easily from this. 
\end{proof}

\begin{defn} \label{defn:basis}
We denote by  $H\stackrel{\sigma}{\hookrightarrow} \Gamma(E\otimes\omega(*D))$
be the $K$-vector subspace 
with basis
\begin{gather}\label{3.3}\frac{e_\mu dt}{(t-a_i)^r},\ \begin{cases}2\le
r\le m_i & 1\le i\le N-1 \\
2\le r\le m_N-1 & i=N \end{cases}, \\
 e_\mu
dt\Big(\frac{1}{t-a_i}-\frac{1}{t-a_N}\Big),\ i< N \notag.
\end{gather}
\end{defn}
There are two splittings of $\sigma$, 
\begin{gather*}\pi_\gamma : \Gamma(E\otimes\omega(*D)) \to
\Gamma(E\otimes\omega(*D))/ \gamma_K\Gamma(E(*D)) \cong H \\
\pi_\nabla : \Gamma(E\otimes\omega(*D)) \to
\Gamma(E\otimes\omega(*D))/ \nabla_{/K}\Gamma(E(*D)) \cong H
\end{gather*}

There is a multiplication map 
$$\eta:\Gamma(E\otimes\omega(*D))\to
\Gamma(E\otimes\omega(*D))\otimes\Omega^1_K.$$ 

We define now two $K$-linear operators.
\begin{defn}\label{defn:op}
The composite map
$$\eta_\gamma := \pi_\gamma\circ \eta \circ \sigma : H
\to H\otimes\Omega^1_K$$
(resp.
$$ \eta_\nabla := \pi_\nabla\circ \eta \circ
\sigma : H \to H\otimes\Omega^1_K)$$
will be called the {\it Higgs} (resp. {\it de Rham}) operator.
Similarly, if $h$ is one of the terms
$\frac{\eta_r^{(i)}}{(t-a_i)^r}$ appearing in the definition of
$\eta$, we denote by $h_\gamma$ and $h_\nabla$ the corresponding
Higgs and de Rham operators.
\end{defn}

The rest of this section is devoted to the comparison of the
trace of those two $K$-linear operators. We will show
\begin{thm}\label{thm:trop}
Let $(E,\nabla)$ be a pseudo-admissible connection on the trivial
bundle $E\cong \oplus_1^r \sO_{\P^1_K}$ on $\P^1_K$ 
having  at least one singularity of order $\ge 2$.
Then
$$\tr (\eta_\gamma -  \eta_\nabla) \equiv  \frac{1}{2} \sum_{m_i\ge 2} m_i
d\log(\det(g^{(i)}_{m_i})) \mod d\log K^{\times}.$$
\end{thm}

We assume henceforth that $m_N\ge 2$. 

Suppose first $a_i$ is a pseudo-logarithmic point for the connection. 
We write
$h=\frac{\eta^{(i)}}{t-a_i}$, and we compute
$\tr(h_\nabla) -\tr(h_\gamma)$. 

The notation $x=y+(H)$ will mean $x$ and $y$ differ by an element in $H$. The
pattern is then we take $x$ in the basis of $H$. We write
\begin{gather}h(x) = \gamma_K y +(H)=\nabla_{/K}y'+(H).
\end{gather}
Then 
\begin{gather}z:=h_\nabla(x) -h_\gamma(x)= -dy'+\gamma_K(y-y').
\end{gather}
Of course, if
$h(x)\in H$ then 
$y=y'=z=0$. Also, we are only interested in the trace, so if the
expansion of $z$ in the basis of $H$ does not involve $x$, we can
ignore it. Suppose e.g. $x=\frac{e_\mu dt}{(t-a_j)^r},\ j\neq i$. Then
$h(x) = \frac{(*)dt}{(t-a_i)(t-a_j)^r}\in H$ (the condition to lie in
$H$ amounts to a bound on the pole order together with no pole of the
differential form at infinity.) Thus such elements $x$ contribute
$0$. Similarly, if $j\neq i$ then 
\begin{gather}\frac{(*)dt}{(t-a_i)(t-a_j)} -
\frac{(*)dt}{(t-a_i)(t-a_N)} \in H
\end{gather} so 
\begin{gather}x=e_\mu dt\Big(\frac{1}{t-a_j}
- \frac{1}{t-a_N}\Big)
\end{gather} 
contributes $0$.    

It remains to consider $x=e_\mu dt\Big(\frac{1}{t-a_i}
- \frac{1}{t-a_N}\Big)$. We have
\begin{gather}h(x) = \frac{\eta_1^{(i)}(e_\mu)dt}{(t-a_i)^2} + (H).
\end{gather}
So we can take
\begin{gather}
y = \frac{(g^{(i)}_1)^{-1}\eta^{(i)}_1(e_\mu)dt}{t-a_i} ;\quad y' =
\frac{(g^{(i)}_1-I)^{-1}\eta^{(i)}_1(e_\mu)dt}{t-a_i} .
\end{gather}
Write $\gamma_K = \frac{g^{(i)}_1dt}{t-a_i} + \gamma_K'$. Then
\begin{gather}z = \gamma_K'(y-y').
\end{gather}
Since we are interested in the trace, we need only consider the
coefficient of $e_\mu dt\Big(\frac{1}{t-a_i}
- \frac{1}{t-a_N}\Big)$ in the expansion of $z$ in the basis of
$H$. This coefficient is the coefficient of $e_\mu$ in 
\begin{gather}\gamma_K'|_{t=a_i}(y-y') = \sum_{\substack{j\neq i\\
r}}(a_i-a_j)^{-r}g^{(j)}_r\Big((g^{(i)}_1 - 
I)^{-1}- (g^{(i)}_1)^{-1}\Big) \eta^{(i)}_1(e_\mu).
\end{gather}
Summing over $\mu$ yields finally
\begin{gather}\label{hhhdr}
\tr (h_\nabla -h_\gamma)=\tr\Big(\sum_{\substack{j\neq i\\
r}}(a_i-a_j)^{-r}g^{(j)}_r\Big((g^{(i)}_1 - 
I)^{-1}- (g^{(i)}_1)^{-1}\Big) \eta^{(i)}_1\Big).
\end{gather}

Next we consider $i$ with $m_i\ge 2$. Take first
$$x = \frac{e_\mu dt}{(t-a_i)^r};\ \begin{cases}2\le
r\le m_i & 1\le i\le N-1 \\
2\le r\le m_N-1 & i=N \end{cases}.
$$
Since we have already handled the $h=\frac{\eta^{(i)}}{t-a_i}$
in pseudo-logarithmic points, 
we introduce the notation
\begin{gather}
\eta'=\eta- \sum_{{\rm pseudo-log}} \frac{\eta^{(i)}}{t-a_i},\\
\gamma'_K=\gamma_K - \sum_{{\rm pseudo-log}}
\frac{g_i^{(i)}}{(t-a_i)}dt. \notag
\end{gather}

Take \begin{equation}\label{3.5} y_0=(g^{(i)}_{m_i})^{-1}
\eta^{(i)}_{m_i-1}(e_\mu)/(t-a_i)^{r-1}.
\end{equation}
It follows from lemma \ref{lem3.3}
that $\eta' x = \gamma'_K y_0 + \text{ lower order terms}$. Here ``lower
order terms'' means terms with denominators $(t-a_j)^s$ where $s\le
m_j,\ j\neq i$ and $s\le m_i+r-2$ for $j=i$. Now continue in this way,
replacing $y_0$ by 
\begin{equation}\label{3.6} y=y_0+\sum_{s=0}^{r-2} \frac{v_s}{(t-a_i)^s}.
\end{equation} 
We may write
\begin{equation}\label{3.7} \eta'(e_\mu)dt/(t-a_i)^r = \gamma'_K y +
\sum_{j=1}^{N-1}\sum_{u=1}^{m_j} \frac{w_{j,u}dt}{(t-a_j)^u} +
\sum_{u=1}^{m_N-1}\frac{w_{N,u}dt}{(t-a_N)^u}.  
\end{equation}
{F}rom equations \eqref{3.5} and \eqref{3.6} we may also write
\begin{multline}\label{3.8} \frac{\eta'(e_\mu)dt}{(t-a_i)^r} = 
(d+\gamma'_K) y +
\sum_{j=1}^{N-1}\sum_{u=1}^{m_j} \frac{w_{j,u}dt}{(t-a_j)^u} +
\sum_{u=1}^{m_N-1}\frac{w_{N,u}dt}{(t-a_N)^u} \\
+ (r-1)(g^{(i)}_{m_i})^{-1}\eta^{(i)}_{m_i-1}(e_\mu)dt /(t-a)^r +
\sum_{s=1}^{r-2} \frac{sv_s dt}{(t-a_i)^{s+1}}.  
\end{multline}

For $r=1$ we may write for suitable $v,w_{j,u} \in \Gamma(E)$
\begin{equation}\label{3.9} \eta'(e_\mu)dt(\frac{1}{t-a_i}-\frac{1}{t-a_N}) =
\gamma'_Kv + \sum_{j=1}^{N-1}\sum_{u=1}^{m_j} \frac{w_{j,u}dt}{(t-a_j)^u} +
\sum_{u=1}^{m_N-1}\frac{w_{N,u}dt}{(t-a_N)^u}.
\end{equation}
Since $dv=0$ in
\eqref{3.9} we conclude from equations
\eqref{hhhdr},
\eqref{3.8} and \eqref{3.9} that 
\begin{multline}\label{3.10} \tr(\eta_\nabla - \eta_\gamma) =
\sum_{i=1}^{N-1} \sum_{r=2}^{m_i}
(r-1)\tr \Big((g^{(i)}_{m_i})^{-1}\eta^{(i)}_{m_i-1}\Big) \\
+  \sum_{r=2}^{m_N-1}
(r-1)\tr \Big((g^{(N)}_{m_N})^{-1}\eta^{(N)}_{m_N-1}\Big)\\
+ \sum_{\substack{i \\ m_i=1}}\tr\Big(\sum_{\substack{j\neq i\\
r}}(a_i-a_j)^{-r}g^{(j)}_r\Big((g^{(i)}_1 - 
I)^{-1}- (g^{(i)}_1)^{-1}\Big) \eta^{(i)}_1\Big).
\end{multline}
(Notice that replacing $\gamma_K$ by $\gamma_K'$ in \eqref{3.7}, \eqref{3.8}
and \eqref{3.9} will not affect the trace calculation.) 
We now use the verticality  condition $dA=A\wedge A \ {\rm mod 
\  } \Omega^2_K\otimes K(X)$.  
\begin{multline} 0 = \sum_{i=1}^N \sum_{r=1}^{m_i}
\frac{dg^{(i)}_r dt}{(t-a_i)^r} + \sum_{i=1}^N \sum_{s=1}^{M_i}
\frac{s\eta^{(i)}_s dt}{(t-a_i)^{s+1}} \\
+ \sum_{i,j=1}^N
\sum_{r,s=1}^{m_i,M_j}
[g^{(i)}_r,\eta^{(j)}_s]\frac{dt}{(t-a_i)^r(t-a_j)^s}+ \sum_{i=1}^N 
\sum_{r=1}^{m_i}[g^{(i)}_r,\eta_0]\frac{dt}{(t-a_i)^r}   .
\end{multline}

Dropping terms with poles at $t=a_i$ of degree $> M_i+1$, we find
\begin{multline} 0 = \sum_{i=1}^N \sum_{r=1}^{m_i}
\frac{dg^{(i)}_r dt}{(t-a_i)^r} + \sum_{i=1}^N \sum_{s=1}^{M_i}
\frac{s\eta^{(i)}_s dt}{(t-a_i)^{s+1}} \\
+ \sum_{i,j=1;\ i\neq j}^N
\sum_{r,s=1}^{m_i,M_j} 
[g^{(i)}_r,\eta^{(j)}_s]\frac{dt}{(t-a_i)^r(t-a_j)^s} \\
 + \sum_{i=1}^N
\sum_{r+s\le M_i+1}[g^{(i)}_r,\eta^{(i)}_s]\frac{dt}{(t-a_i)^{r+s}} 
+  \sum_{i=1}^N  \sum_{r=1}^{m_i}[g^{(i)}_r,\eta_0]\frac{dt}{(t-a_i)^r} .
\end{multline}

In an admissible point, we calculate as before, multiplying through by
$$(g^{(i)}_{m_i})^{-1}(t-a_i)^{m_i}/dt$$ and set $t=a_i$ to 
get
\begin{multline}\label{eq3.27} 0 = (g^{(i)}_{m_i})^{-1} dg^{(i)}_{m_i} +
(m_i-1)(g^{(i)}_{m_i})^{-1}\eta^{(i)}_{m_i-1} \\
+ \sum_{j\neq i}
\sum_{s=1}^{M_j} (g^{(i)}_{m_i})^{-1}[g^{(i)}_{m_i},\eta^{(j)}_s]/(a_i-a_j)^s
+ \sum_{r+s=m_i}(g^{(i)}_{m_i})^{-1}[g^{(i)}_r,\eta^{(i)}_s] \\
+ (g^{(i)}_{m_i})^{-1}[g^{(i)}_{m_i},\eta_0] .
\end{multline}

Taking traces gives
\begin{multline}\label{3.28} 0 = d\log(\det(g^{(i)}_{m_i})) +
(m_i-1)\tr((g^{(i)}_{m_i})^{-1}\eta^{(i)}_{m_i-1}) \\
+\sum_{r+s=m_i}\tr\Big((g^{(i)}_{m_i})^{-1}[g^{(i)}_r,\eta^{(i)}_s]\Big) .
\end{multline}

On the other hand,  in a pseudo-logarithmic point,
we multiply through by $(t-a_i)^2/dt$
and then set $t=a_i$ getting 
\begin{gather}\label{3.29}\eta^{(i)}_1 = [\eta^{(i)}_1,g^{(i)}_1].
\end{gather}
Now discard those terms, multiply by $(t-a_i)/dt$ and set $t=a_i$. One gets
\begin{gather}\label{3.30} 0=dg^{(i)}_1 + \sum_{j\neq i}\sum_{s=1}^{M_j}
(a_i-a_j)^{-s}\Big([g^{(i)}_1,\eta_s^{(j)}] +
[g^{(j)}_s,\eta^{(i)}_1]\Big) + [g^{(i)}_1,\eta_0].
\end{gather}
Formula \eqref{3.29} gives
\begin{gather}
\eta^{(i)}_1(g^{(i)}_1-I)=g^{(i)}_1\eta^{(i)}_1,
\end{gather}
whence, assuming the indicated matrices invertible, one has
\begin{gather}
\eta^{(i)}_1(g^{(i)}_1-I)^{-1} =(g^{(i)}_1)^{-1}\eta^{(i)}_1.
\end{gather}
Using $\tr(a[b,c])=0$ if $[a,b]=0$, multiplying equation \eqref{3.30} on
the left by $(g^{(i)}_1)^{-1}$ (resp. by $(g^{(i)}_1-I)^{-1}$) and
taking traces yields
\begin{gather}\tr\Big(\sum_{j\neq i}\sum_{s=1}^{M_j}
(a_i-a_j)^{-s}(g^{(i)}_1)^{-1}
[g^{(j)}_s,\eta^{(i)}_1]\Big) \in d\log K^\times, \\
\tr\Big(\sum_{j\neq i}\sum_{s=1}^{M_j} 
(a_i-a_j)^{-s}(g^{(i)}_1-I)^{-1}
[g^{(j)}_s,\eta^{(i)}_1]\Big) \in d\log K^\times.\notag
\end{gather} 

We now get
\begin{multline}\label{3.34} \sum_{\substack{i \\ m_i=1}}\tr\Big(
\sum_{\substack{j\neq i\\ r}} (a_i-a_j)^{-r}g^{(j)}_r\Big((g^{(i)}_1 -   
I)^{-1}- (g^{(i)}_1)^{-1}\Big) \eta^{(i)}_1\Big) \\
= \sum_{\substack{i \\ m_i=1}}\tr\Big( \sum_{\substack{j\neq i\\
r}} (a_i-a_j)^{-r}\Big((g^{(i)}_1 -   
I)^{-1}- (g^{(i)}_1)^{-1}\Big) \eta^{(i)}_1g^{(j)}_r\Big) \\ 
\equiv \sum_{\substack{i \\ m_i=1}}\tr\Bigg( \sum_{\substack{j\neq i\\
r}} (a_i-a_j)^{-r}\Big(\eta^{(i)}_1(g^{(i)}_1 -   
I)^{-1}- (g^{(i)}_1)^{-1}
\eta^{(i)}_1\Big)g^{(j)}_r\Bigg) \\
\equiv 0 \mod(d\log K^\times). 
\end{multline}

Now one can compare \eqref{3.10}, \eqref{3.28}, and \eqref{3.34} and deduce:

\begin{prop}\label{prop:hidr}
With notation as in definition \ref{defn:op} above, and assuming that $a_N$ has
multiplicty $m_N\ge 2$, one has
\begin{multline*}\tr(\eta_{\nabla} - \eta_\gamma) \equiv
-\sum_{i=1}^{N-1}\frac{m_i}{2}\Big(d\log(\det(g_{m_i}^{(i)}))
+ \tr\sum_{r+s=m_i} (g^{(i)}_{m_i})^{-1}[g^{(i)}_r,\eta^{(i)}_s]\Big) \\
- \frac{m_N-2}{2} \Big(d\log(\det(g_{m_N}^{(N)}))+ \tr
\sum_{r+s=m_N} (g^{(N)}_{m_N})^{-1}[g^{(N)}_r,\eta^{(N)}_s]\Big) \\
\mod d\log K^\times.  
\end{multline*}
\end{prop} 

To complete the proof of theorem \ref{thm:trop}, we must show
\begin{prop}\label{prop:iden}
Let $(E,\nabla)$ be a vertical admissible connection, with local equation
$$\frac{g_m dz}{z^m} + \frac{g_{m-1} dz}{z^{m-1}}+ \ldots + 
\frac{\eta_{m-1}}{z^{m-1}}+ \frac{\eta_{m-2}}{z^{m-2}} + \ldots $$
in an admissible point. Then
$$\Phi=\tr(g_m^{-1}\sum_{s=0}^{m-1}[g_{m-s},\eta_s]) = 0
$$
\end{prop}
\begin{proof}
The connection $\nabla$ is vertical. Vanishing for curvature
terms involving $z^{-p}$ for $p\ge m+1$ implies
\begin{gather}\label{3.15}
[g_m, \eta_{\ell}] + [g_{m-1}, \eta_{\ell + 1}] + \ldots
[g_{\ell + 1}, \eta_{m-1}]=0 \\
{\rm for \ } \ell=1, \ldots, m-1 .\notag
\end{gather}

The tactic is to eliminate first $\eta_1$ from $\Phi$, then
$\eta_2$ etc. One easily verifies matrix relations
\begin{gather}\label{3.17}
[a^{-1}, b]= -a^{-1}[a,b]a^{-1}\\
\tr (g^{-1}[a,b])=-\tr (a[g^{-1},b])=\tr (ag^{-1}[g,b]g^{-1}) \notag \\
\tr (ag^{-1}[b,\eta]g^{-1})=\tr (ag^{-1}b\eta g^{-1}) -\tr
(bg^{-1}ag^{-1}\eta). \notag \end{gather}
In particular, 
\begin{gather}
\tr (ag^{-1}[b,\eta]g^{-1})=\tr (ag^{-1}bg^{-1}[g,\eta]g^{-1})\ 
{\rm if \ }\  ag^{-1}b=bg^{-1}a.
\end{gather}
Write 
\begin{gather}
\Phi= \tr g_m^{-1}[g_{m-1}, \eta_1] + \tr
\sum_{s=2}^{m-1}g_m^{-1}[g_{m-s}, \eta_s].
\end{gather}
This yields
\begin{gather}
\Phi= -\tr g_m^{-1}g_{m-1}g_{m}^{-1}[g_{m-1},\eta_2] + \tr
g_m^{-1}[g_{m-2}, \eta_2] \\
- \tr \sum_{s=3}^{m-1}
g_m^{-1}g_{m-1}g_m^{-1}[g_{m-s+1}, \eta_s] +
\tr \sum_{s=3}^{m-1}g_m^{-1}[g_{m-s}, \eta_s].\notag
\end{gather} 
Applying again now the relations \eqref{3.17} yields
\begin{gather}
\Phi= \tr(-g_m^{-1}g_{m-1}g_m^{-1}g_{m-1}g_m^{-1} +
g_m^{-1}g_{m-2}g_m^{-1})[g_m, \eta_2] \\
- \tr \sum_{s=3}^{m-1}
g_m^{-1}g_{m-1}g_m^{-1}[g_{m-s+1}, \eta_s] +
\tr \sum_{s=3}^{m-1}g_m^{-1}[g_{m-s}, \eta_s].\notag
\end{gather}
Assume inductively that for some $t\ge 2$, one can write $\Phi$ as follows:
\begin{gather}\label{3.22}
\Phi=\tr  (\sum_{a=1}^t (-1)^{a-1} \sum_{\tau_1+ \ldots +
\tau_a=t} g_m^{-1}g_{m-\tau_1} \cdots g_m^{-1}g_{m-\tau_a}g_m^{-1})
[g_m, \eta_t] \\ + \tr
\sum_{s=t+1}^{m-1} \sum_{\ell=0}^{t-1} (\sum_{a=0}^{\ell} (-1)^a
\sum_{\tau_1 + \ldots +
\tau_a=\ell}g_m^{-1}g_{m-\tau_1} \cdots
g_m^{-1}g_{m-\tau_a}g_m^{-1}) [g_{m-s+\ell},\eta_s].\notag 
\end{gather}
Applying \eqref{3.15} to the first line, and isolating the terms in
$\eta_{t+1}$ and in $\eta_s, s\ge (t+2)$, one obtains
\begin{gather}
\Phi= F(t+1) +\\
\tr
\sum_{s=t+2}^{m-1} \sum_{\ell=0}^{t} (\sum_{a=0}^{\ell} (-1)^a
\sum_{\tau_1 + \ldots +
\tau_a=\ell}g_m^{-1}g_{m-\tau_1} \cdots
g_m^{-1}g_{m-\tau_a}g_m^{-1}) [g_{m-s+\ell},\eta_s],\notag 
\end{gather}
with
\begin{gather}
F(t+1)=
\tr (\sum_{a=1}^t (-1)^a \sum_{\tau_1+\ldots \tau_a=t}
g_m^{-1}g_{m-\tau_1} \cdots
g_m^{-1}g_{m-\tau_a}g_m^{-1})[g_{m-1}, \eta_{t+1}] \\
+ \tr \sum_{\ell=0}^{t-1}(\sum_{a=0}^\ell
(-1)^a\sum_{\tau_1+\ldots +\tau_a=\ell}
g_m^{-1}g_{m-\tau_1} \cdots
g_m^{-1}g_{m-\tau_a}g_m^{-1})[g_{m-t-1+\ell},\eta_{t+1}]. \notag
\end{gather}
It remains to arrange $F(t+1)$. To this aim, write
\begin{gather}
F(t+1)= \sum_{\ell=0}^t \sum_{a=0}^\ell \sum_{\tau_1 +\ldots
+\tau_a=\ell} \tr \\
((-1)^a g_{m-\tau_1}g_m^{-1}\cdots
g_{m-t-1+\ell}\eta_{t+1}g_m^{-1} - (-1)^a g_{m-t-1+\ell}
g_m^{-1}\cdots g_{m-\tau_a}g_m^{-1}\eta_{t+1}).\notag
\end{gather}
Now we group those terms differently. To a tuple 
$(\tau_1,\ldots,\tau_a)$, with $\tau_1+\ldots +
\tau_a=\ell$, we associate the tuple 
$(\tau'_1,\ldots,\tau'_{a})$ with $\tau'_1+\ldots \tau'_{a} +
\tau_1=t+1$, $\tau'_{a}=t+1-\ell$, and otherwise
$\tau_i=\tau'_{i-1}$ for $i\ge 2$. Using the first relation of
\eqref{3.15} again, this gives for those 2 terms  together
\begin{gather}
\tr (-1)^{a} g_m^{-1}g_{m-\tau_1}g_m^{-1}\cdots
g_{m-t-1+\ell}g_m^{-1}[g_m,\eta_{t+1}].
\end{gather}
This shows that the relation \eqref{3.22} is true, with $t$ replaced by
$t+1$. As the last equation of \eqref{3.15} for $\ell= m-1$ is $[g_m,
\eta_{m-1}]=0$,  
one obtains by induction that $\Phi$ vanishes on the variety defined
by \eqref{3.15}.
\end{proof}

\begin{remark}\label{rmk3.4} Writing $a_i=g_i$ and $b_i = \eta_{i-1}$,
the above 
proposition can be restated as follows. Suppose $a(t) =
a_mt^m+\ldots+a_1t$ and $b(t) = b_mt^m+\ldots+b_1t$ are polynomials
with matrix coefficients satisfying $a(0)=b(0)=0$ and $a_m$
invertible. Assume $[a(t),b(t)]=c_mt^m + \text{ lower order
terms}$. Then $\tr(a_m^{-1}c_m) = 0$. 
\end{remark}

\section{The Gau\ss-Manin Determinant: Step 1}\label{sect4}

In this section we begin the computation of the Gau\ss-Manin
determinant appearing in the main theorem \ref{thm:mainthm}.

 We keep the same notations as in section 3. 
In particular, $E$ is a trivial bundle on $\P^1_K$ with
basis $e_\mu$, having at least one point of multiplicity $\ge 2$, 
$D=\{a_1,\dotsc,a_N\}$, and $H\hookrightarrow
\Gamma(E\otimes \omega(*D))$ is 
the $K$-subspace with basis 
defined in \ref{defn:basis}.
We continue to write $\nabla
= d+\gamma+\eta$ and $\nabla_{/K}=d+\gamma_K$ as in definition
\ref{defn:gammaK}.

The Gau\ss-Manin connection is computed from the diagram

\begin{equation}\label{4.1} \minCDarrowwidth.5cm \begin{CD}
@. \Gamma(E(*D)) @= \Gamma(E(*D)) \\ 
 @. @VV\nabla V @VV\nabla_{/K}V \\
 \Gamma(E(*D))\otimes \Omega^1_K @>>> \Gamma(E\otimes\Omega^1(*D))
@>\stackrel{s}{\leftarrow} >> \Gamma(E\otimes \omega(*D)) \\
 @VV\nabla_{/K}\otimes 1V  @VV\nabla V \\
 \Gamma(E\otimes\omega(*D))\otimes \Omega^1_K  @>\cong
>\iota >\Gamma(E\otimes \Omega^2(*D)/F^2)  
\end{CD}
\end{equation}

Here in the central column $\Omega$ refers to the K\"ahler
differentials $\Omega_{\P^1_K/k}$, and $F^2 := \sO_{\P^1}(*D)\otimes
\Omega^2_K \subset \Omega^2(*D)$. We are interested in the induced map
from $H^1_{DR}=H^1_{DR/K}(\P^1-D, (E,\nabla))$, which is the cokernel
of the right hand column, to $H^1_{DR}\otimes_K \Omega^1_K$, which is
the cokernel of the left hand column. Let $p:\Gamma(E\otimes
\omega(*D)) \to H^1_{DR}$ be the projection. By construction, $H
\hookrightarrow \Gamma(E\otimes \omega(*D))$ splits $p$. Let
$q:H^1_{DR} \cong H$ denote the splitting. 
The section $s$ is given on $H$ by
\begin{multline}\label{4.2} s\Big(\frac{e_\mu dt}{(t-a_i)^r}\Big) =
\frac{e_\mu d(t-a_i)}{(t-a_i)^r}; \\
s\Big(e_\mu dt\big(\frac{1}{t-a_i}-\frac{1}{t-a_N} \big)\Big) = \frac{e_\mu
d(t-a_i)}{t-a_i}-\frac{e_\mu
d(t-a_N)}{t-a_N}.
\end{multline}

\begin{lem}
The Gau\ss-Manin determinant is the trace of the map 
$$(q\otimes 1)\circ (p\otimes 1)\circ \iota^{-1}\circ\nabla\circ s : H
\to H.
$$
\end{lem}
\begin{proof}Straightforward. \end{proof}
Explicitly, this map is obtained by applying the projection $(q\otimes
1)\circ (p\otimes 1)$ to the right hand side in
\begin{multline}\label{4.3} \frac{e_\mu dt}{(t-a_i)^r}\mapsto
\sum_{j=1}^N\sum_{s=1}^{m_j} \frac{g^{(j)}_s(e_\mu) d(a_i-a_j)\wedge
dt}{(t-a_j)^s(t-a_i)^r} \\
+\sum_j\sum_{s=1}^{M_j}\frac{\eta^{(j)}_s(e_\mu)\wedge
dt}{(t-a_j)^s(t-a_i)^r }  
\end{multline}
\begin{multline}\label{4.4}  e_\mu dt\Big(\frac{1}{t-a_i}-\frac{1}{t-a_N}\Big)
\mapsto \\
\sum_{j=1}^N\sum_{s=1}^{m_j}g^{(j)}_s(e_\mu) \Big
( \frac{d(a_i-a_j)\wedge dt}{(t-a_j)^s(t-a_i)}- \frac{d(a_N-a_j)\wedge
dt}{(t-a_j)^s(t-a_N)}\Big) \\
+   \sum_j\sum_{s=1}^{M_j}\eta^{(j)}_s(e_\mu)\wedge
dt\Big(\frac{1}{(t-a_j)^s(t-a_i)}- \frac{1}{(t-a_j)^s(t-a_N)}\Big).  
\end{multline}

We leave aside for the moment the terms in the trace involving $\eta$
and focus on the trace of the map which we rewrite using lemma
\ref{lem3.3} as
\begin{multline}\label{4.5} \frac{e_\mu dt}{(t-a_i)^r} \mapsto \sum_{j=1}^N
\sum_{s=1}^{m_j} g^{(j)}_s(e_\mu)d(a_i-a_j)\wedge
dt\Big(\frac{(a_i-a_j)^{-s}}{(t-a_i)^r}+\ldots \Big) 
\end{multline}
\begin{multline}\label{4.6} e_\mu
dt\Big(\frac{1}{t-a_i}-\frac{1}{t-a_N}\Big) \mapsto \\
\sum_{j=1}^N \sum_{s=1}^{m_j}
g^{(j)}_s(e_\mu)\Bigg[(a_i-a_j)^{-s}d(a_i-a_j)\wedge dt\Big 
 ( \frac{1}{t-a_i}-\frac{1}{t-a_N}\Big) \\  
+ \Big((a_N-a_j)^{-s}d(a_N-a_j)\wedge dt -
 (a_i-a_j)^{-s}d(a_i-a_j)\wedge dt\Big) \\
\times \Big(\frac{1}{t-a_j}-
\frac{1}{t-a_N}\Big)+\ldots  \Bigg].
\end{multline}
Terms are to be dropped if some factor becomes $0$.
The terms represented by ellipses $(\ldots)$ do not enter into the
trace calculation. Also all terms lie in the $K$-span of the basis
\eqref{3.5}. Let $\Psi$ denote the resulting endomorphism of $H$. The
contribution to $\tr(\Psi)$ from \eqref{4.5} is 
\begin{multline}\sum_{i=1}^N (m_i-1)\sum_{\substack{j=1 \\  j\neq
i}}^N\sum_{s=1}^{m_j} 
\tr(g^{(j)}_s)(a_i-a_j)^{-s}d(a_i-a_j)\wedge dt \\
- \sum_{j=1}^{N-1}\sum_{s=1}^{m_j}\tr(g^{(j)}_s)(a_N-a_j)^{-s}d(a_N-a_j)
\wedge dt.    
\end{multline}
The second sum arises because $\frac{e_\mu dt}{(t-a_N)^{m_N}}$ is not
a basis element for $H$. The contribution from \eqref{4.6} is
\begin{multline}\sum_{i=1}^{N-1}\sum_{\substack{j=1 \\ j\neq i}}^N
\sum_{s=1}^{m_j} \tr(g^{(j)}_s)(a_i-a_j)^{-s}d(a_i-a_j)\wedge dt \\
+ \sum_{j=1}^{N-1} \sum_{s=1}^{m_j}
\tr(g^{(j)}_s)(a_N-a_j)^{-s}d(a_N-a_j) \wedge dt. 
\end{multline}
In total, this gives
\begin{multline}\label{4.7} \sum_{i=1}^N m_i\sum_{\substack{j=1 \\
j\neq i}}^N \sum_{s=1}^{m_j} 
\tr(g^{(j)}_s)(a_i-a_j)^{-s}d(a_i-a_j)\wedge dt \\
- \sum_{j=1}^{N-1}
\sum_{s=1}^{m_j}\tr(g^{(j)}_s)(a_N-a_j)^{-s}d(a_N-a_j) 
\wedge dt.  
\end{multline}

In fact, the above analysis of $\tr(\Psi)$ omits some terms. On the
right in \eqref{4.3} taking $j=N$ and $s=m_N$ gives a term
$$\frac{g^{(N)}_{m_N}(e_\mu)d(a_i-a_N)\wedge dt}{(t-a_N)^{m_N}(t-a_i)^r}.
$$  
Expanding this by lemma \ref{lem3.3} yields a term (for $2\le r\le m_i$)
\begin{multline}\label{4.8}
\frac{g^{(N)}_{m_N}(e_\mu)(a_N-a_i)^{-r}d(a_i-a_N)\wedge 
dt}{ (t-a_N)^{m_N}} \\
\equiv -\sideset{}{'}\sum_{j,s}
\frac{g^{(j)}_s(e_\mu)(a_N-a_i)^{-r}d(a_i-a_N)\wedge dt}{ (t-a_j)^s}.
\end{multline}
Here $\equiv$ means equivalent in $H^1_{DR}\otimes\Omega^1_K$. The
prime in the sum 
means omit the pair $j=N,\ s=m_N$.  

Similarly, from \eqref{4.4} we get a term
\begin{multline}\label{4.9} \frac{g^{(N)}_{m_N}(e_\mu)(a_N-a_i)^{-1}
d(a_i-a_N)\wedge dt}{(t-a_N)^{m_n}} \\
\equiv  -\sideset{}{'}\sum_{j,s}
\frac{g^{(j)}_s(e_\mu)(a_N-a_i)^{-1}d(a_i-a_N)\wedge dt}{ (t-a_j)^s}. 
\end{multline} 

Of course, in \eqref{4.8} and \eqref{4.9} the contribution to the trace
comes from $j=i$. These precisely cancel the second double sum on the
right in \eqref{4.9}. Thus, one gets
\begin{equation}\label{4.12} \tr(\Psi) = \sum_{i=1}^N m_i\sum_{\substack{j=1 \\
j\neq i}}^N \sum_{s=1}^{m_j} 
\tr(g^{(j)}_s)(a_i-a_j)^{-s}d(a_i-a_j)\wedge dt. 
\end{equation}

Finally, comparing \eqref{4.3}, \eqref{4.4}, and definition
\ref{defn:op}, we have
\begin{prop}\label{prop4.1} With notation as above (definition
\ref{defn:op} and equation
\eqref{4.12}), 
the Gau\ss-Manin trace on $H^1_{DR}$ is given by
\begin{multline*}\tr_{GM}(H^1_{DR}) = \tr(\Psi) + \tr(\eta_\nabla) \\
\stackrel{{\rm thm.}\ref{thm:trop}}{\equiv}
\tr(\Psi) + \tr(\eta_\gamma)+\frac{1}{2}\sum_{i;\ m_i\ge 2}m_i
d\log(\det(g^{(i)}_{m_i})) \mod d\log K^\times. 
\end{multline*}
\end{prop}

\section{The Higgs Trace}\label{sect5}

The purpose of this section and of the next one 
is to rewrite the Higgs trace
$\tr(\eta_\gamma)$ as a sum of terms which are in some sense local,
associated to the singularities $a_i$ of the connection. In the next
section we will compute these local traces. 

It will be convenient to write
\begin{gather}\label{5.1} \sV := \Gamma(E(*D));\quad  \sW :=
\Gamma(E\otimes \omega(*D))\subset \sV dt ;\\ 
\gamma_K = \sum_{j=1}^N
\sum_{r=1}^{m_j}\frac{g_r^{(j)}dt}{(t-a_j)^r}:\sV \to \sW. \notag
\end{gather}
(See equations \eqref{strict} and definition \ref{defn:gammaK}).

We identify $H\hookrightarrow \sW$ with basis \eqref{3.3}. As in
section \ref{sect3}, there is a splitting $H\cong \sW/\gamma_K\sV$. We
write
\begin{equation}\label{5.2}\eta = \eta_0 + \sum_{j=1}^N
\sum_{s=1}^{M_j}\frac{\eta^{(j)}_s}{(t-a_j)^s} = \eta_0 + \sum_j
\eta^{(j)}. 
\end{equation}
The $\eta^{(j)}_s$ are matrices with entries in $\Omega^1_K$. We view
these objects as linear maps $H \to H\otimes_K \Omega^1_K$:
$$H \hookrightarrow \sW \stackrel{\eta^{(i)}}{\longrightarrow} \sW
\otimes \Omega^1_K  
\twoheadrightarrow H\otimes \Omega^1_K
$$

Now fix an $i$. It will be convenient to put $a_i$ at $\infty$, 
so we set
\begin{gather}\label{coord}
u:= \frac{1}{t-a_i}.
\end{gather} 
Define
\begin{multline}\label{5.3} R:=\Gamma(\P^1-D,\sO) =
K[\frac{1}{t-a_1},\dotsc,\frac{1}{t-a_N}] \\
= K\Big[\frac{u}{1-(a_1-a_i)u},\dotsc,u,\dotsc,\frac{u}{1-(a_N-a_i)u}\Big].
\end{multline}
In the $u$ coordinates, 
\begin{equation}\label{5.4}\gamma_K = \sum
g_s^{(j)}\Big(\frac{u}{1-(a_j-a_i)u}\Big)^s.  
\end{equation}
The basis of $H$ is
\begin{gather}\label{5.5} e_\mu dt\Big(\frac{u}{1-(a_j-a_i)u}\Big)^r;
\quad 2\le r\le m_j \text{ (resp. $2\le r\le m_N-1$)}, \\  
e_\mu dt\Big(\frac{u}{1-(a_j-a_i)u} -
\frac{u}{1-(a_N-a_i)u}\Big) \notag \\
=e_\mu
dt\frac{(a_j-a_N)u^2}{(1-(a_j-a_i)u)(1-(a_N-a_i)u)} ;\quad j\neq N. \notag
\end{gather}
Define 
\begin{gather}\label{5.6} \theta := \prod_{j\neq i}
(1-(a_j-a_i)u)^{m_j}, \\
\theta_1 := \theta/(1-(a_N-a_i)u), \notag \\
g = \theta \cdot \gamma_K = \sum_{j,s} g_s^{(j)}u^s. \notag
\end{gather}
Note that $\theta$ is a unit in $R$. We can write
\begin{multline}\label{5.7} g = \prod_{j\neq i}
(a_i-a_j)^{m_j}g^{(i)}_{m_i} u^m + \text{lower order terms} \\
= u^2+\text{higher order terms}, 
\end{multline}
where $m=\sum m_j$. Let $V=\oplus_\mu Ke_\mu = \Gamma(\P^1,E)$ and
write  
\begin{equation}\label{5.8} H' := Vu^2 dt \oplus\ldots \oplus Vu^{m-1} dt
\subset V[u]dt.
\end{equation}
As a consequence of \eqref{5.8} we have
\begin{equation}\label{5.9} \theta_1 H = H'.
\end{equation}
We are interested in the trace of $\eta^{(i)} =
\sum_{s=1}^{M_i}\eta^{(i)}_s u^s$. 
Consider the diagram
\begin{equation}\label{5.10} \minCDarrowwidth.3cm\begin{CD}H @>\subset >>
\sW @>\eta^{(i)}>> \sW \otimes \Omega^1_K @>>> \sW /\gamma_K\sV \otimes
\Omega^1_K @<\cong << H\otimes \Omega^1_K \\
@VV\cong V @V\theta_1 V\cong V @A\theta_1^{-1}A\cong A @A\theta_1^{-1}
A\cong A @AA \cong A  \\ 
H' @>\subset>> \sW @>\eta^{(i)}>> \sW \otimes \Omega^1_K
@>>> \sW /\gamma_K\sV \otimes
\Omega^1_K @<\cong << H'\otimes \Omega^1_K \\
@| @A\cup AA @A\cup AA @A\cong A a A @| \\
H' @>\subset>> u^2V[u] dt @>\eta^{(i)}>> u^2V[u] dt\otimes \Omega^1_K @>>>
\frac{u^2V[u]}{gV[u]}dt\otimes \Omega^1_K  @<\cong << H'\otimes \Omega^1_K
\end{CD}
\end{equation}

\begin{prop}\label{prop5.1}The following maps have the same trace
\begin{gather*} \eta^{(i)} : H \to H\otimes \Omega^1_K \\
\eta^{(i)} : H' \to H'\otimes \Omega^1_K \\
\eta^{(i)} : \frac{V[u]}{gV[u]}dt \to \frac{V[u]}{gV[u]}dt\otimes \Omega^1_K.
\end{gather*}
Here the first two maps are given by horizontal rows in
\eqref{5.10}. The third is given by embedding
$$\frac{V[u]}{gV[u]}dt \hookrightarrow V[u]dt
$$
via the basis $V\oplus Vu \oplus\ldots\oplus Vu^{m-1}$ and then
proceeding as in \eqref{5.10}.  
\end{prop}
\begin{proof} The first two traces are equal by the diagram. For the
third, note that since $g=u^2+ \text{ higher}$, it follows that one
has an exact sequence compatible with the endomorphism multiplication
by $u$
$$0 \to H' \to V[u]/gV[u]dt \to  V[u]/u^2V[u]dt \to 0.
$$
Since $\eta^{(i)}$ has no constant term in $u$, it acts nilpotently on
the right.
\end{proof}

\section{The Higgs trace: local calculation}\label{sect6}

In this section we give a formula for the trace of $\eta^{(i)}$ as in
proposition \ref{prop5.1}, involving residues. As already mentioned in the 
introduction, the method here is reminiscent of the
classical residue calculation for the trace of an element in a field
extension. 

To simplify, we write $h$ in place of
$\eta^{(i)}$. We also suppress the $\Omega^1_K$ and treat $h(u)$ as a
polynomial with matrix coefficients. The case of matrices with
coefficients in $\Omega^1_K$ follows immediately by applying arbitrary
derivations $\Omega^1_K \to K$ to the entries.

To avoid confusion we
write 
\begin{equation}\label{6.1} \phi(h) : V[u]/gV[u] \hookrightarrow V[u]
\stackrel{h\cdot}{\longrightarrow} V[u] \twoheadrightarrow V[u]/gV[u], 
\end{equation}
where $ V[u]/gV[u] \hookrightarrow V[u]$ is defined as in
proposition \ref{prop5.1} via the invertibility of the leading
coefficient of $g$
\begin{gather}
V\oplus Vu \oplus\ldots\oplus Vu^{m-1}\hookrightarrow V[u].
\end{gather}

By \eqref{4.7} and admissibility, the leading coefficient of $g(u)$ is
invertible so $g^{-1}\in {\rm End}(V)((u^{-1}))$. Write
$dg=\frac{dg}{du}du$ and let 
$$\res\!\!_{u=\infty} :  {\rm End}(V)((u^{-1})) \to
{\rm End}(V)$$ 
be the evident extension of the residue map.

One has
\begin{prop}\label{prop6.1} The notations being  as above, one has
$$ \tr_{V[u]/gV[u]}(\phi(h)) = - \tr_V \res\!\!_{u=\infty}(dgg^{-1}h).
$$
\end{prop}
\begin{proof}Write $g=a_0u^m + a_1 u^{m-1} +\ldots + a_m$. Note that
neither side of the identity changes if we replace $g$ by $ga_0^{-1}$
so we may assume $a_0=1$. 

Also, by linearity, we may assume
$h=cu^p$. The matrix
for the action of $u$ on $V[u]/gV[u]$, the entries of which are themselves
matrices, is
$$M = \begin{pmatrix}0 & 0 & \hdots & -a_n \\
1 & 0 & \hdots  & -a_{n-1} \\
\vdots & \vdots & \vdots & \vdots\\
0 & 0 & 1 & -a_1
\end{pmatrix}
$$
The matrix for $u^p$ is $M^p$. We write $\tr M^p\in {\rm End}(V)$ 
for the naive
trace, i.e. the sum of the diagonal elements. E.g. $\tr M = -a_1\in
{\rm End}(V)$. Then 
$$\tr(\phi(cu^p)) = \tr_V(\tr(M^p)c).
$$ 
Also
$$\tr_V \res\!\!_{u=\infty}(dgg^{-1}cu^p) = \tr_V( \res(dgg^{-1}u^p)\cdot c).
$$
(The residue is computed in the ring ${\rm End}(V)((u^{-1}))$. Since $u$ is in
the center of this ring, we may move $c$ past $u^p$ under the residue.) It
will therefore suffice to show
\begin{equation}\label{6.2} \tr(M^p) = -\res\!\!_{u=\infty}(dgg^{-1}u^p).
\end{equation}
Let $z=u^{-1}$ and write $g=u^mG(z)$ with $G(z) = I+a_1z+\ldots+a_mz^m$.
The assertion becomes
\begin{equation}\label{6.3} dGG^{-1} = (\tr(M)+ \tr(M^2)z+\ldots)dz.
\end{equation}

\begin{lem}\label{lem6.2} Let $X = (x_{ij})$ be an $m\times m$ matrix. Then
\begin{equation}\tr(XM^p) = \sum_{q=0}^p(-1)^q
\sum_{\substack{1\le m_1,\dotsc,m_q\le m\\ 1\le i,i_1\le m \\ 
\sum m_k=p+i-i_1 \\
m_1\ge m-i_1+1}}x_{i,i_1}a_{m_1}\ldots a_{m_q}.
\end{equation}
In particular, taking $X=I$ it follows that
\begin{equation}\label{6.5} \tr(M^p) = \sum_{q=1}^p(-1)^q
\sum_{\substack{1\le m_1,\dotsc,m_q\le m\\ 1\le i\le m \\ 
\sum m_k=p \\
m_1\ge m-i+1}}a_{m_1}\ldots a_{m_q}.
\end{equation}
\end{lem} 
\begin{proof}[proof of lemma] Write $M= (M_{ij})_{1\le i,j\le m}$. We
have $M_{i+1,i}=1,\ M_{i,m}=-a_{m+1-i}$, and $M_{ij} = 0$ otherwise.
Thus
\begin{multline} \tr(XM^p) = \sum_{i,i_1,\dotsc,i_p}
x_{i,i_1}M_{i_1,i_2}\cdots M_{i_p,i} = \sum_{i,i_1,q,j_1,\dotsc,j_q}
x_{i,i_1}M_{j_1,m}\cdots M_{j_q,m} \\
= \sum_{i,i_1,q,j_1,\dotsc,j_q}(-1)^q x_{i,i_1}a_{m-j_1+1}\cdots
a_{m-j_q+1}
\end{multline}
The conditions on the tuples $\{i,i_1,q, j_1,\dotsc,j_q\}$ over which the
right hand sums are taken become
\begin{gather}\label{6.7} 0\le q\le p;\ j_1\le i_1;\\
(i_1-j_1)+1+(m-j_2)+1+\ldots+(m-j_q)+1+(m-i) = p \notag
\end{gather}
Replacing $j_k$ by $m_k:=m-j_k+1$, these become
\begin{equation}\label{6.8}0\le q\le p;\ m_1\ge m-i_1+1;\ \sum m_k = p+i-i_1,
\end{equation}
proving the lemma. 
\end{proof}

Write $c_p = \tr(M^p)$. We must show
\begin{multline}\label{6.9}-(c_1+c_2z+c_3z^2+\ldots)(1+a_1z+\ldots+a_mz^m) \\
= a_1+2a_2z+\ldots + ma_mz^{m-1}. 
\end{multline}
This amounts to
\begin{equation}\label{6.10}c_p+c_{p-1}a_1+\ldots+c_1a_{p-1} =
\begin{cases}-pa_p & p<m \\
0 & \text{else}\end{cases} 
\end{equation}
Suppose first $p<m$. With reference to \eqref{6.5}, one can isolate the
terms in $c_p$ ending in $a_k$ for $1\le k\le p-1$ and write
\begin{equation}c_p = c_{p,1}a_1+c_{p,2}a_2+\ldots+c_{p,p-1}a_{p-1} + R_p.
\end{equation}
Here 
\begin{equation}c_{p,k} =  -\sum_{r=1}^{p-k}(-1)^r \sum_{\substack{1\le
m_1,\dotsc,m_r\le m\\ 1\le i\le m \\ 
\sum m_k=p-k \\
m_1\ge m-i+1}}a_{m_1}\ldots a_{m_r} = -c_{p-k}.
\end{equation}
(Notice that since each $m_j\ge 1$, terms with $r>p-k$ are impossible.
Also, since $k<p$, necessarily $r\ge 1$.) The remainder $R_p$ is given by
the terms $a_p$ in $c_p$. In the sum for $c_p$ these terms arise when
$q=1$ and $p\ge m-i+1$, i.e. $m-p+1\le i\le m$. There are $p$ such terms:
\begin{equation}R_p = - a_p.
\end{equation}
Finally, in \eqref{6.10} we consider the terms with $p\ge m$. Writing $t_k =
\tr(XM^k)$ and replacing $X$ with $M^j$ for some $j$, it suffices to show
\begin{equation}t_0a_m+\ldots + t_{m-1}a_1+t_m = 0.
\end{equation}
We start with
\begin{equation}t_m = \sum_{q=1}^p(-1)^q
\sum_{\substack{1\le m_1,\dotsc,m_q\le m\\ 1\le i,i_1\le m \\ 
\sum m_k=n+i-i_1 \\
m_1\ge n-i_1+1}}x_{i,i_1}a_{m_1}\ldots a_{m_q}.
\end{equation}
Note $q=0$ is not possible because $\sum m_k = m+i-i_1\ge 1$. Again, by
grouping together the terms ending with $a_k$ we get
$$ t_m = -t_{m-1}a_1-\ldots -t_0 a_m,
$$
which is the desired equation. 
\end{proof}

\begin{cor}\label{cor6.3} We have with notation as in definition \ref{defn:op}
\begin{multline*}\tr(\eta_{\nabla}) \equiv (m-2)\tr(\eta_0)-\sum_i
\res_{t=a_i}\tr (dgg^{-1}\eta^{(i)}) \\
+ \frac{1}{2}\sum_{i;\ m_i\ge 2} d\log(\det(g^{(i)}_{m_i}))
\mod d\log(K^\times).  
\end{multline*}
\end{cor}
\begin{proof}By theorem \ref{thm:trop}, we may replace $\eta_\nabla$
with $\eta_\gamma$. By proposition \ref{prop4.1} $\tr(\eta_\gamma-\eta_0)$ is
the sum of the $\tr(\eta^{(i)})$ on $V[u]/gV[u]dt$. By proposition
\ref{prop5.1} this is the same as $-\sum_i\tr_V
\res_{t=a_i}(dgg^{-1}\eta^{(i)})$. Note the factor $m-2$ on the right
is because as a matrix $\eta_0$ acts on $V=\Gamma(\P^1,E)$ while the
trace one wants is the action on $H \cong V^{\oplus m-2}$.  
\end{proof}

\section{The proof of the main Theorem}\label{sect7}

In this section we deduce the main theorem \ref{thm:mainthm}
from the equality of the Higgs and de Rham
traces (theorem \ref{thm:trop}) and from the shape of the Higgs trace 
(proposition \ref{prop6.1}) via residues.

We start with an admissible connection $(E, \nabla)$ on $\P^1_K$.
Let $V \subset U$ be a Zariski open subset, with complement $Z$. By
localization, one obtains
\begin{gather}
\det H_{DR}(U, \nabla_{/K})=\det H_{DR}(V, \nabla_{/K}) +
\det \nabla|_{Z}.
\end{gather}
On the other hand,   at a special pseudo-logarithmic point the local
factor $\tr(dgg^{-1}A)=0$. Indeed, 
writing the connection as $g_1\frac{dz}{z} +
g_0dz + \ldots + \frac{\eta_1}{z} + \frac{\eta_0} + \ldots$, the
local factor is $\tr g_0g_1^{-1} \eta_1$. 
The special pseudo-logarithmic points
have, in the notations of the definition \ref{defn:spe},
 a local matrix of the shape
\begin{gather}
\Big(\begin{array}{ll}
A+m\frac{dz}{z} & zB\\
\frac{C}{z} & D +n\frac{dz}{z}
\end{array}\Big).
\end{gather}
Thus in particular
\begin{gather}
\eta_1= 
\Big(\begin{array}{ll}
0 & 0\\
\gamma_0 & 0
\end{array}\Big),
\end{gather}
where $C=cdz + \gamma_0 + \gamma_1 z + \ldots$,
while $g_0$ and $g_1$ are both of the shape
\begin{gather}
\Big(\begin{array}{ll}
* & 0\\
* & *
\end{array}\Big).
\end{gather}
Thus $\tr g_0g_1^{-1}\eta_1=0$.

Thus the difference of the
right hand side of the theorem \ref{thm:mainthm} for $U$ and $V$ is
the difference of the global factors, which is
$\det \nabla|_{Z}$, as one sees taking a trivializing
section which is good for $V$. Thus by theorem \ref{thm:red}, we may assume
that $E= \oplus_1^r \sO_{\P^1_K}$.

Let $G=\gamma_K\cdot\prod_j
(t-a_j)^{m_j}$. Note $G=u^{-m}g(u)$ with $u$ 
as in \eqref{coord}. Write the absolute connection as $\nabla = d+A$
with $A=\gamma + \eta$ as in definition \ref{defn:gammaK}. 
\begin{prop}\label{prop7.1} We have
\begin{multline}\label{7.1} \sum_i \tr\ \res\!\!_{t=a_i}dGG^{-1}A = 
-\sum_{\substack{i,j,r\\j\not = i}}\tr
(g_r^{(i)})m_j(a_j-a_i)^{-r}d(a_j-a_i) \\
+\sum_i \tr\ \res\!\!_{t=a_i}(dGG^{-1}\eta^{(i)}).
\end{multline}
\end{prop}
\begin{proof}Define absolute
forms $s$, $s(i)$ and $\eta(i)$:
\begin{equation}\label{7.2}s = \frac{dt}{\prod_j
(t-a_j)^{m_j}};\quad s(i) = \frac{d(t-a_i)}{\prod_j
(t-a_j)^{m_j}};\quad \nabla = d +G\cdot  
s(i) +\eta(i).
\end{equation}

The local term at $t=a_i$ is
$$\tr\ \res\!\!_{t=a_i}dGG^{-1}A=\tr\
\res\!\!_{t=a_i}dG\cdot s(i)+\tr\ \res\!\!_{t=a_i}dGG^{-1}\eta(i).
$$
Applying trace to $dA= A\wedge A$ yields
$$\tr(dG\cdot s(i) + G\cdot ds(i) + d\eta(i)) = 0.
$$
(Note here that $s(i)$ is not closed as an absolute form!) Also, modulo
$\Omega^2_K$, we have $\res\ \tr(d\eta(i)) = 0$ because the residue of
an exact form vanishes. The local
term thus becomes
\begin{equation}\label{7.3}\tr\ \res\!\!_{t=a_i}dGG^{-1}A= -\tr\
\res\!\!_{t=a_i}G\cdot ds(i) + \tr\ \res\!\!_{t=a_i}dGG^{-1}\eta(i).
\end{equation}
We have
\begin{multline}G\cdot ds(i) = \\
\Big(\prod_j
(t-a_j)^{m_j}\sum_{r,k}\frac{g_r^{(k)}}{(t-a_k)^r}\Big)\frac{\sum_j
m_j\prod_{k\neq j}(t-a_k)d(t-a_j)\wedge d(t-a_i)}{\prod_j(t-a_j)^{m_j+1}}
\\
= \Big(\sum_{r,k} \frac{g^{(k)}_r}{(t-a_k)^r}\Big)\sum_j m_j
\frac{d(a_i-a_j)\wedge dt}{t-a_j}.
\end{multline}
Thus
\begin{equation}\label{7.5}\res\!\!_{t=a_i}\ \tr(G\cdot ds(i)) =
\res\!\!_{t=a_i}\Big(\sum_{\substack{j,k,r\\ j\not =i}} \frac{\tr(g^{(k)}_r)m_j
d(a_i-a_j)\wedge dt}{(t-a_k)^r(t-a_j)}\Big). 
\end{equation}
Expanding
$$\frac{1}{(t-a_j)} = -(a_j-a_i)^{-1}\sum_{n=0}^\infty
\Big(\frac{t-a_i}{a_j-a_i}\Big)^n
$$
and substituting on the right in \eqref{7.5} 
\begin{equation}\label{7.6}\res\!\!_{t=a_i}\ \tr(G\cdot ds(i)) =
\sum_{\substack{j,r\\j\not = i}}\tr
(g_r^{(i)})m_j(a_j-a_i)^{-r}d(a_j-a_i)
\end{equation}
Thus, \eqref{7.5} becomes after summing over $i$
\begin{multline}\label{7.7} \sum_i \tr\ \res\!\!_{t=a_i}dGG^{-1}A = 
-\sum_{\substack{i,j,r\\j\not = i}}\tr
(g_r^{(i)})m_j(a_j-a_i)^{-r}d(a_j-a_i) \\
+\sum_i \tr\ \res\!\!_{t=a_i}(dGG^{-1}\eta(i)).
\end{multline}

It remains to compare $\eta^{(i)}$ and $\eta(i)$. Recall
$\gamma_K=G\cdot s$. We have
\begin{gather}\gamma = \gamma_K - \sum_{j,r}\frac{g^{(j)}_r
da_j}{(t-a_j)^r} \\
G\cdot s(i) = \gamma_K - \sum_{j,r}\frac{g^{(j)}_r
da_i}{(t-a_j)^r} \\
\nabla = d+G\cdot s(i) + \eta(i) = d+\gamma + \eta_0 + \sum_j \eta^{(j)}.
\end{gather}
{F}rom these equations it follows that
$$\eta(i) - \eta^{(i)} = \sum_{j\neq i}\frac{k^{(j)}_r}{(t-a_j)^r}.
$$
Since this difference has no pole at $a_i$, and since $G(a_i)$ is
invertible, we find 
$$\res\!\!_{t=a_i}(dGG^{-1}\eta(i)) = \res\!\!_{t=a_i}(dGG^{-1}\eta^{(i)})
$$
Making this substitution in \eqref{7.7} proves the proposition. 
\end{proof}

Our calculations to this point have been on $H^1_{DR}$ which
introduces a minus sign in the final formula. (Note, admissibility
forces $H^0_{DR}=(0)$.) We find therefore
\begin{multline}\label{7.11} -\tr_{GM}(H^*_{DR}) =
\tr(\Psi)+\tr(\eta_{\nabla})\quad\text{(proposition \ref{prop4.1})} \\
\equiv \tr(\Psi) + (m-2)\tr(\eta_0) - \sum_i \tr\ 
\res\!\!_{t=a_i}(dgg^{-1}\eta^{(i)}) \\
+ \frac{1}{2}\sum_{i;\ m_i\ge 2} d\log(\det(g^{(i)}_{m_i}))
\mod  d\log K^\times \quad \text{(corollary \ref{cor6.3})} \\
\equiv \tr(\Psi) + (m-2)\tr(\eta_0) - \sum_i \tr\ 
\res\!\!_{t=a_i}(dGG^{-1}\eta^{(i)}) \\
+ \frac{1}{2}\sum_{i;\ m_i\ge 2} d\log(\det(g^{(i)}_{m_i}))
\mod  d\log K^\times \\
\equiv (m-2)\tr(\eta_0) - \sum_i \tr\ 
\res\!\!_{t=a_i}(dGG^{-1}A)\\
+ \frac{1}{2}\sum_{i;\ m_i\ge 2} d\log(\det(g^{(i)}_{m_i}))
\mod  d\log K^\times \quad \text{(proposition \ref{prop7.1})}. 
\end{multline}

Note here we can replace $g$ by $G=u^{-m}g$ because
$\eta^{(i)}$ only involves terms of degrees $\ge 1$ in $u$, so 
$\res(\frac{du}{u}\eta^{(i)}) = 0$.

View $s=\frac{dt}{\prod_j (t-a_j)^{m_j}}$ as a
relative form, i.e. as a meromorphic section of $\omega(\sum
m_j(a_j))$. As such, it has a zero of order $m-2$ at infinity,
$(s)=(m-2)\infty$. Note $\det(E,\nabla)|_\infty = (K,\tr(\eta_0))$. 
Since the relative connection is given by $\nabla_{/K} = G\cdot s$,
the desired formula is 
\begin{multline}\label{7.12} \tr_{GM}(H^*_{DR}) \equiv
-\det(E,\nabla)|_{(s)} + \sum_i  \res\!\!_{t=a_i} \tr (dGG^{-1}A) \\
+ \frac{1}{2}\sum_{i;\ m_i\ge 2} d\log(\det(g^{(i)}_{m_i}))
\mod  d\log K^\times. 
\end{multline}

Finally, comparing \eqref{7.11} and \eqref{7.12} we see that theorem
 \ref{thm:mainthm} 
holds for pseudo-admissible connections on bundles on
$\P^1$.

\bibliographystyle{plain}
\renewcommand\refname{References}

\end{document}